\documentclass[11pt]{amsart}

\usepackage{amsmath,amsthm, amscd, amssymb, amsfonts}
\usepackage{epsfig}
\input xy
\xyoption{all}

\newcommand{\Hc}{{\mathcal H}}

\newcommand{\Mo}{{\mathcal M}}
\newcommand{\No}{{\mathcal N}}
\newcommand{\mo}{{\mathcal M}}

\newcommand{\emod}{{\,}_S^H{\mathcal M}_K}

\newcommand{\Ss}{{\mathcal S}}
\newcommand{\ot}{{\otimes}}

\newcommand{\ere}{{\mathcal R}}
\newcommand{\ele}{{\mathcal L}}
\newcommand{\el}{{\underline L}}
\newcommand{\elt}{{\overline L}}
\newcommand{\ca}{{\mathcal C}}

\newcommand{\action}{\leftharpoonup}

\newcommand{\Fc}{{\mathcal F}}
\newcommand{\Gc}{{\mathcal G}}

\newcommand{\Ker}{\mbox{\rm Ker\,}}
\newcommand{\op}{\rm{op}}
\newcommand{\cop}{\rm{cop}}
\newcommand{\bop}{\rm{bop}}

\newcommand{\ku}{{\Bbbk}}

\newcommand{\Na}{{\mathbb N}}

\newcommand{\uno}{{\mathbf 1}}

\newcommand{\id}{\mbox{\rm id\,}}

\newcommand{\Vect}{\mbox{\rm Vect\,}}

\newcommand{\Fun}{\operatorname{Hom}}

\newcommand\Rep{\operatorname{Rep}}

\newcommand\Irr{\operatorname{Irr}}
\newcommand\rk{\operatorname{rk}}

\newcommand\co{\operatorname{co}}

\newcommand\Hom{\operatorname{Hom}}
\newcommand\uhom{\underline{\Hom}}

\newcommand{\End}{\operatorname{End}}
\newcommand\uend{\underline{\End}}

\newcommand\Cent{\operatorname{Cent}}
\newcommand\St{\operatorname{Stab}}

\renewcommand{\_}[1]{\mbox{$_{\left( #1 \right)}$}}

\theoremstyle{plain}

\numberwithin{equation}{section}

\newtheorem{teo}{Theorem}[section]

\newtheorem{lema}[teo]{Lemma}

\newtheorem{cor}[teo]{Corollary}

\newtheorem{prop}[teo]{Proposition}

\newtheorem{step}{Step}

\theoremstyle{definition}

\newtheorem{defi}[teo]{Definition}

  \newtheorem{exa}[teo]{Example}

\theoremstyle{remark}

\newtheorem{rmk}[teo]{Remark}

\def\pf{\begin{proof}}

\def\epf{\end{proof}}

\theoremstyle{remark}

\begin{document}

\title[Module categories over  Hopf algebras]
{On module categories over finite-dimensional \\ Hopf algebras}
\author[Andruskiewitsch and Mombelli]{Nicol\'as Andruskiewitsch and
Juan Mart\'\i n Mombelli }
\thanks{This work was partially supported by
Agencia C\'ordoba Ciencia, ANPCyT-Foncyt, CONICET,
TWAS (Trieste), Fundaci\'on Antorchas and Secyt (UNC)}
\address{Facultad de Matem\'atica, Astronom\'\i a y F\'\i sica, Universidad Nacional de C\'ordoba
\newline
\indent CIEM -- CONICET
\newline \indent Medina Allende s/n, (5000) Ciudad Universitaria, C\'ordoba, Argentina}
\email{andrus@mate.uncor.edu, \quad \emph{URL:}\/
http://www.mate.uncor.edu/andrus} \email{mombelli@mate.uncor.edu,
\quad \emph{URL:}\/ http://www.mate.uncor.edu/mombelli}
\begin{abstract} We show that indecomposable exact module
categories over the category $\Rep H$ of representations of a
finite-dimensional Hopf algebra $H$ are classified by left
comodule algebras, $H$-simple from the right and with trivial
coinvariants, up to equivariant Morita equivalence. Specifically,
any indecomposable exact module categories is equivalent to the
category of finite-dimensional modules over a left comodule
algebra. This is an alternative approach to the results of Etingof
and Ostrik. For this, we study the stabilizer introduced by Yan
and Zhu and show that it coincides with the internal Hom. We also
describe the correspondence of module categories between $\Rep H$
and $\Rep (H^*)$.
\end{abstract}

\subjclass{16W30, 18D10, 19D23}
\date{\today}
\maketitle

\setcounter{tocdepth}{1} \tableofcontents

\section*{Introduction}

The notion of fusion category is a far-reaching generalization of
the notion of finite group. It has been studied in connection with
different problems in conformal field theory, mechanical
statistics, the theory of subfactors, and others, the common theme
being ``quantum symmetries". A comprehensive presentation of
fusion categories is \cite{ENO}, see also references therein. We
refer to \emph{loc. cit.} for definitions and notations used in
the present paper. There is a notion of ``module category over a
tensor category" known in category theory since the sixties
\cite{Be}.  Semisimple module categories over a fusion category
should play the same fundamental r\^ole as the representation
theory of finite groups; see \cite{O1, O2}. In the beautiful paper
\cite{eo}, the notion of \emph{finite tensor category} was
introduced and several properties of fusion categories were
extended to finite tensor categories. Finite tensor categories are
like fusion categories but without semisimplicity; a basic example
is the category of finite-dimensional representations of a
finite-dimensional Hopf algebra. The natural class of module
categories over a finite tensor category is the class of exact
tensor categories \cite{eo}.
 Let $\ca$ be a finite tensor category. Then

    \smallbreak
\begin{itemize}
    \item Any exact module category over $\ca$
    is a finite direct product of  \emph{indecomposable} exact module
    categories.

    \smallbreak

    \item An indecomposable module category $\Mo$ over $\ca$
    is naturally equivalent-- as a module category-- to
    the category of right $A$-modules in $\ca$, where $A$ is an
    algebra in $\ca$ \cite{O2, eo}. Explicitly, $A$ can be
    chosen as the internal End of any non-zero object in
    $\Mo$, \emph{cf.} the proof of \cite[Thm. 1]{O2}.

\end{itemize}

\smallbreak The purpose of this paper is to study exact module
categories over the finite tensor category $\ca = \Rep H$, where
$H$ is a finite-dimensional Hopf algebra over algebraically closed
field $\ku$ of characteristic 0. In \cite{eo}, the authors propose
to classify indecomposable exact module categories over $\ca$ by
classifying first simple from the right exact $H$-module algebras.
Indeed, any indecomposable exact module category over $\Rep H$ is
equivalent to the category $_H\Mo_R$ of Hopf $(H, R)$-bimodules
for some $H$-module algebra $R$. Instead, our approach is through
left comodule algebras: any indecomposable exact module category
over $\Rep H$ is also equivalent to the category $_K\Mo$ of left
modules over some $H$-comodule algebra $K$. We feel that this
approach is intuitively clearer. Both approaches are related
because of the correspondence between module categories over $\Rep
H$ and $\Rep (H^*)$.

Section \ref{section:mc} contains a general discussion of module
categories arising from comodule algebras. We show that
indecomposable exact module category are classified by exact
indecomposable left $H$-comodule algebras up to a suitable
``equivariant" Morita equivalence, see Theorem \ref{clasif-uno}.
Examples of indecomposable left $H$-comodule algebras are the
coideal subalgebras, thanks to a recent result of Skryabin, see
Proposition \ref{Hsimplealgexact}.

Section \ref{section:yz} is devoted to the stabilizer introduced
by Yan and Zhu and a generalization thereof. We study this
construction and prove a general version of the duality, Theorem
\ref{cor:yzduality-simple}, answering a question of Yan and Zhu
\cite[p. 3897]{YZ}. We also generalize a formula for the dimension
of the stabilizer obtained by Zhu for semisimple Hopf algebras
\cite{Z}.

In Section 3, we show that the Yan-Zhu stabilizer is the internal
Hom and apply our version of the Yan-Zhu duality to prove a finer
classification of indecomposable exact module categories in terms
of comodule algebras, see Theorem \ref{mod=lca}. We also describe
explicitly the correspondence between module categories over $\Rep
H$ and $\Rep (H^*)$ using the Yan-Zhu stabilizer, see Theorem
\ref{corr4}.

\subsection*{Acknowledgements} We thank Sonia Natale for her
interest in this work and many important remarks, and
Hans-J\"urgen Schneider for help with the proof of Prop.
\ref{Hsimplealgexact}.

\section{Module categories}\label{section:mc}

\subsection{Preliminaries}
Let $\ku$ be an algebraically closed field of characteristic 0.
All vector spaces, algebras, unadorned $\otimes$ and Hom, are over
$\ku$. As usual, $V^* = \Hom(V, \ku)$ is the dual vector space of
a vector space $V$, and $\langle\,,\,\rangle: V^* \times V \to
\ku$ is the evaluation. If $V$ is finite dimensional, we identify
$(V\otimes V)^*$ with $V^*\otimes V^*$ via
$$
\langle\alpha \otimes \beta, v\otimes w\rangle = \langle\alpha,
v\rangle \langle\beta, w\rangle,
$$
$\alpha, \beta\in V^*$, $v,w\in V$. If $A$ is an algebra, then
${}_A\Mo$, resp. $\Mo_A$ denotes the category of
finite-dimensional left, resp. right, $A$-modules. If $C$ is a
coalgebra, then ${}^C\Mo$, resp. $\Mo^C$ denotes the category of
finite-dimensional left, resp. right, $C$-comodules. The kernel of
the counit of $C$ is denoted $C^+$. We shall use Sweedler's
notation for the coproduct and coactions: $\Delta(c) = c\_1 \ot
c\_2$ if $c\in C$. Also, if $\lambda: M \to C\ot M$ is a left
coaction, $\lambda(m) = m\_{-1} \ot m\_0$.

Let $H$ be a finite-dimensional Hopf algebra with multiplication
$\mu$, comultiplication $\Delta$, counit  $\varepsilon$, antipode
$\Ss$. We denote $\Rep H$ instead of ${}_H\Mo$ to emphasize the
presence of the tensor structure. We denote by $L: H\to \End H$
and by $R:H\to \End H$  the left regular representation, resp. the
righ regular representation, that is $L_a(b) = ab$, $R_a(b)=ba$
for $a,b\in H$. Recall that $H^{\op}$, resp. $H^{\cop}$ is the
Hopf algebra with opposite multiplication, resp. opposite
comultiplication, and $H^{\bop} = (H^{\op})^{\cop}$. Clearly,
$(H^*)^{\cop} = (H^{\op})^*$.

We denote by $\rightharpoonup: H\ot H^*\to H^*$ and
$\leftharpoonup: H^*\ot H\to H^*$ the actions obtained by
transposition of the right and left multiplications, and by
$\leftharpoondown: H^*\ot H\to H^*$ and $\rightharpoondown: H\ot
H^*\to H^*$ the corresponding compositions with the inverse of the
antipode. That is,

\begin{align} \label{arponalto}
\langle a\rightharpoonup\alpha,b\rangle &= \langle
\alpha,ba\rangle = \langle \alpha\leftharpoonup b,a\rangle,
\\ \label{arponbajo} \langle b\rightharpoondown\alpha, a\rangle
&= \langle \alpha, \Ss^{-1} (b)a \rangle, \qquad \langle
\alpha\leftharpoondown a, b\rangle = \langle \alpha, b\Ss^{-1}(a)
\rangle,
\end{align}
$a,b \in H$, $\alpha\in H^*$. We denote by $\el: H\to \End (H^*)$,
respectively by $\elt: H\to \End (H^*)$, the representation
afforded by $\rightharpoonup$, respectively by
$\rightharpoondown$. The analogous actions (and representations)
of $H^*$ on $H$ are denoted by the same symbols. With respect to
$\rightharpoonup$, $H^*$ is a $H$-module algebra. Note that
\begin{align*}
\alpha \rightharpoonup h &= \langle \alpha, h_{(2)}\rangle \,
h_{(1)}, &\quad  h\leftharpoonup \alpha &= \langle \alpha,
h_{(1)}\rangle \,
h_{(2)},\\
\alpha\rightharpoondown h&=\langle\alpha,\Ss^{-1}(h\_1) \rangle
\,h\_2, &\quad h\leftharpoondown \alpha&= \langle\alpha,
\Ss^{-1}(h\_2)\rangle \,h_{(1)},
\end{align*}
$\alpha\in H^*$, $h\in H$. Notice that
\begin{align} \label{yz-tecnico3}(ht)\leftharpoondown \alpha
= (h\leftharpoondown \alpha\_2)(t\leftharpoondown \alpha\_1).
\end{align}

\smallbreak  If $X, Y\in {}_H\Mo$ then $\Hom(X, Y)$ is again a
left $H$-module via
\begin{align} \label{accionhomizq}
(h\cdot T)(x) = h_{(1)} \cdot T (\Ss(h_{(2)}) \cdot x), \qquad
x\in X, \, h\in H, \, T\in \Hom(X, Y).
\end{align}

Similarly,  if $W, Z\in \Mo_H$ then $\Hom(W, Z)$ is again a right
$H$-module via
\begin{align} \label{accionhomder}
(T\cdot h) (w) =  T(w \cdot \Ss^{-1}(h_{(2)}))\cdot h_{(1)},
\qquad w\in W, \, h\in H, \, T\in \Hom(W, Z).
\end{align}

\begin{lema}\label{lema:yz-tecnico} \cite{YZ}. If $h, t\in H$, $\alpha, \beta\in H^*$ then
\begin{align} \label{yz-tecnico}
(h \rightharpoondown \beta)\alpha &= (h \leftharpoonup
\Ss^{-2}(\alpha_{(1)})) \rightharpoondown (\beta\alpha_{(2)}).
\qed\end{align}
\end{lema}

\smallbreak Let $H$ be a Hopf algebra. A (left) $H$-comodule
algebra is a (left) comodule that is also an associative unital
algebra such that the coaction and the counit are morphisms of
algebras. A (left) $H$-module algebra is a (left) module $A$ that
is also an associative unital algebra such that the coaction and
the counit are morphisms of algebras. This means that $h\cdot(ab)
= (h\_1\cdot a)(h\_2\cdot b)$ and $h\cdot 1= \varepsilon(h)1$, for
all $h\in H$, $a,b\in A$. Right module algebras are defined
similarly. If $V\in {}_H\Mo$ (resp. $V\in \Mo_H$) then the algebra
$\End(V)$ is a left (resp. right) $H$-module algebra. If $H$ is
finite-dimensional, the notions of ``left $H$-module algebra" and
``right $H^*$-comodule algebra" are equivalent.

\smallbreak Let $R$ be an $H$-module algebra. An $(H,R)$-module is
an $R$-module $M$ inside the monoidal category $_{H}\Mo$;  in
other words, $M$ is a left $H$-module and a right $R$-module, and
the right $R$-action $M\ot R\to M$ is an $H$-module map. We denote
by $_{H}\Mo_{R}$ the category of $(H,R)$-modules with morphisms
the maps preserving both actions.

\smallbreak Let $K$ be a left $H$-comodule algebra. In the same
vein as before, we denote by ${}^{H}_{K}\Mo$ the category of left
$H$-comodules, left $K$-modules $M$ such that the left $K$-module
structure $K\ot M\to M$ is an $H$-comodule map. Analogously, if
$S$ is a right $H$-comodule algebra, then there is a category
$_{S}\Mo{}^{H}$ of right $H$-comodules, left $S$-modules with
action being a morphism of $H$-comodules.

\smallbreak Let $K$ be a left $H$-comodule algebra. A left
$H$-ideal $I$ of $K$ is a left ideal that is also an $H$-comodule.
Right and two-sided $H$-ideals are defined similarly. The
following natural notions are discussed in \cite{eo, S}.

\begin{defi}\label{hsimple}  We shall say that $K$ is {\it $H$-simple from the left}, (resp.
{\it $H$-simple from the right}, resp. {\it $H$-simple}) if $K$
has no non-trivial left (resp. right, resp. two-sided) $H$-ideal.
We shall say that $K$ is {\it $H$-indecomposable} if there are no
nontrivial two-sided $H$-ideals $I$ and $J$ such that $K = I
\oplus J$.
\end{defi}

We shall need later the following result.

\begin{lema}\label{inyectiva}
Assume that $K$ is $H$-simple. If $W\neq 0$ is a left $K$-module,
then the representation
\begin{equation}\label{e1}\rho:K\to \End(H^{*}\ot W), \quad \rho(k)(\alpha\ot w)=
k\_{-1}\rightharpoondown \alpha\ot k\_0\cdot  w,\end{equation} is
faithful. Similar for module algebras.
\end{lema}

\begin{proof}  Let $I$ be the two-sided ideal $\Ker(\rho)$. We shall prove that $I$ is a
$H$-ideal. First, notice that if $k\in I$ then
\begin{equation}\label{id} <\gamma,\Ss^{-1}(k\_{-1})>\, k\_0\cdot w=0,
\end{equation} for all
$\gamma\in H^{*}$, $w\in W$.

Let $k\in I$. Then $\lambda(k)\in H\ot I$ iff
$(\id\ot\rho)\lambda(k)=0$ iff $(\Ss^{-1}\ot\rho)\lambda(k)=0$.
Let $\alpha, \beta\in H^{*}, w\in W$. Evaluating
$(\Ss^{-1}\ot\rho)\lambda(k)$ in $(\alpha\ot w)$ and the applying
$\beta\ot \id$ we get
\begin{align*} <\beta,\Ss^{-1}(k\_{-1})>\rho(k\_0)(\alpha\ot w)&=
<\beta,\Ss^{-1}(k\_{-2})><\alpha\_1, \Ss^{-1}(k\_{-1}) >
\alpha\_2\ot k\_0 \cdot w\\
&=<\Ss^{-1}(\alpha\_1\beta),k\_{-1}>\,\alpha\_2\ot k\_0 \cdot w.
\end{align*}
Evaluating the last expression in $h\ot \id$, $h\in H$ we obtain
\begin{align*} <\Ss^{-1}(\alpha\_1\beta),k\_{-1}><\alpha\_2,h> k\_0 \cdot
w&= <\Ss^{-1}((h\rightharpoonup\alpha)\beta ),k\_{-1}> k\_0 \cdot
w\\
&=<(h\rightharpoonup\alpha)\beta,\Ss^{-1}(k\_{-1})> k\_0 \cdot w.
\end{align*}
The last expression is 0 by \eqref{id}. Since $\alpha, \beta, w$
and $h$ are arbitrary, $(\Ss^{-1}\ot\rho)\lambda(k)=0$.
\end{proof}

\subsection{Freeness over comodule
algebras}\label{subsection:skryabin} Let $H$ be a finite
dimensional Hopf algebra. In this subsection we recall some
important results obtained recently by S. Skryabin \cite{S}.

\begin{teo}\label{Skryabin1}
Let $A$ be a finite dimensional right $H$-comodule algebra. Assume
that $A$ is $H$-simple. Then all objects in $\Mo^H_A$ are
projective $A$-modules.  $M \in \Mo^H_A$ is a free $A$-module if
and only if $M/MQ$ is a free $A/Q$-module for at least one
two-sided maximal ideal $Q$ of $A$.
\end{teo}

\pf This is a particular case of \cite[Theorem 4.2]{S}. \epf

Clearly, similar statements hold also for other combinations like
${}_K^H\Mo$.

\begin{teo}\label{freeness} Let $A$ be a finite dimensional $H$-simple right $H$-comodule algebra.
Let  $M \in \Mo^H_A$. Then there exists $t\in \Na$ such that $M^t$
is a free $A$-module.
\end{teo}

\pf This follows from the proof of \cite[Theorem 3.5]{S}.\epf

\begin{prop}\label{coideal-simple} If $K\subseteq H$ is a left
coideal subalgebra then $K$ is $H$-simple.
\end{prop}
\pf See the proof of \cite[Theorem 6.1]{S}.\epf

\subsection{Tensor and module categories}
In this paper, we stick to the following terminology. We refer to
\cite{BK, O2} for more details.

\begin{itemize}
    \item A \emph{monoidal category } is a category $\ca$ provided with a
    ``tensor" functor $\ot: \ca \times \ca\to \ca$, a unit object
    $\uno \in \ca$, associativity constraint $a$, left and right
    unit constraints $l$ and $r$, all these subject to the
    pentagon and the triangle axioms.

\smallbreak\item A monoidal category $\ca$ is \emph{rigid} if any
object in $\ca$ has right and left duals.

\smallbreak\item A \emph{tensor category } is a rigid monoidal
category with $\ca$ abelian and $\ot$ additive in each variable
(that is, $\ot$ is a bifunctor). We shall assume hereafter that
both $\ca$ and $\ot$ are $\ku$-linear.

\smallbreak\item The  \emph{opposite monoidal category }
$\ca^{\op}$ to a monoidal category $\ca$
  is the same category $\ca$ but with $X\ot^{\op}Y = Y\ot X$,
  $a^{\op}_{X,Y,Z} = a_{Z,Y,X}^{-1}$, $l^{\op} = r$, $r^{\op} =
  l$, and the same unit. If $\ca$ is rigid, resp. tensor, then so
  is $\ca^{\op}$.

\smallbreak\item A \emph{monoidal functor} between monoidal
categories $\ca$ and $\ca'$ is a functor $F:\ca \to \ca'$,
provided with a natural isomorphism $ b_{X,Y}: F(X\otimes Y)\to
F(X)\otimes F(Y)$ and an isomorphism $u: F(\uno) \to \uno$
satisfying two natural axioms, namely
\begin{align*}
a_{F(X),F(Y),F(Z)} (b_{X,Y}\otimes \id)b_{X\otimes Y,Z} &=
(\id\otimes b_{Y,Z})b_{X,Y\otimes Z}F(a_{X,Y,Z})
\\ F(l_X) =  l_{F(X)}(u\otimes \id)b_{\uno,X}, \qquad
F(r_X) &=  r_{F(X)}(\id\otimes u)b_{X,\uno}.
\end{align*}

\smallbreak\item A \emph{tensor functor} between tensor categories
$\ca$ and $\ca'$ is a monoidal functor $F:\ca \to \ca'$ which is
$\ku$-linear.

\smallbreak\item  \cite{Be} A {\em module category} over a tensor
category $\ca$ is an abelian category $\Mo$ provided with an exact
bifunctor\footnote{that is, exact in each variable} $\otimes: \ca
\times \Mo \to \Mo$ and natural associativity and unit
isomorphisms $m_{X,Y,M}: (X\otimes Y)\otimes M \to X\otimes
(Y\otimes M)$, $\ell_M: \uno \otimes M\to M$ such that for any
$X,Y,Z\in \ca$, $M\in \Mo$.

\begin{align}\label{action1}
(\id\otimes m_{Y,Z,M})m_{X,Y\otimes Z,M}(a_{X,Y,Z}\otimes id) &=
m_{X,Y,Z\otimes M}\, m_{X\otimes Y,Z,M}
\\\label{action2}
(\id\otimes \ell_M)m_{X,\uno ,Y} &=r_X\otimes \id
\end{align}

\smallbreak\item   A {\em module functor} between module
categories $\Mo$ and $\Mo'$ over a tensor category $\ca$ is a pair
$(\Fc,c)$, where $\Fc:\Mo \to \Mo'$ is a $\ku$-linear functor and
$c_{X,M}: \Fc(X\otimes M)\to X\otimes \Fc(M)$ is a natural
isomorphism such that the following diagrams are commutative, for
any $X, Y\in \ca$, $M\in \Mo$:

\begin{align}\label{modfunctor1}
(\id_X\otimes c_{Y,M})c_{X,Y\otimes M}\Fc(m_{X,Y,M}) &=
m_{X,Y,\Fc(M)}\, c_{X\otimes Y,M}
\\\label{modfunctor2}
\ell_{\Fc(M)} \,c_{\uno ,M} &=\Fc(\ell_{M}).
\end{align}
We shall denote $(\Fc, c): \Mo \to \Mo'$. Note that $(\id, \id):
\Mo \to \Mo$ is a module functor. There is a composition of module
functors: if $\Mo''$ is another module category and $(\Gc,d): \Mo'
\to \Mo''$ is another module functor then the composition
\begin{equation}\label{modfunctor-comp}
(\Gc\circ \Fc, e): \Mo \to \Mo'', \qquad \text {where } e = d\circ
\Gc(c),
\end{equation} is
also a module functor.

\smallbreak\item  Let $\Mo_1$ and $\Mo_2$ be module categories
over $\ca$. We denote by $\Hom_{\ca}(\Mo_1, \Mo_2)$ the category
whose objects are module functors $(\Fc, c)$ from $\Mo_1$ to
$\Mo_2$. A morphism between  $(\Fc,c)$ and
$(\Gc,d)\in\Hom_{\ca}(\Mo_1, \Mo_2)$ is a natural transformation
$\alpha: \Fc \to \Gc$ such that for any $X\in \ca$, $M\in \Mo_1$:
\begin{gather}
\label{modfunctor3} d_{X,M}\alpha_{X\otimes M} =
(\id_{X}\ot\alpha_{M})c_{X,M}.
\end{gather}

\smallbreak\item  Two module categories $\Mo_1$ and $\Mo_2$ over
$\ca$ are {\em equivalent} if there exist module functors
$F:\Mo_1\to \Mo_2$ and $G:\Mo_2\to \Mo_1$ and natural isomorphisms
$\id_{\Mo_1} \to F\circ G$, $\id_{\Mo_2} \to G\circ F$ that
satisfy \eqref{modfunctor3}.

\smallbreak\item  The {\em direct sum} of two module categories
$\Mo_1$ and $\Mo_2$ over a tensor category $\ca$  is the
$\ku$-linear category $\Mo_1\times \Mo_2$ with coordinate-wise
module structure.

\smallbreak\item  A module category is {\em indecomposable} if it
is not equivalent to a direct sum of two non trivial module
categories.
\end{itemize}

Let $\Irr \Mo$ be the set of isomorphism classes of irreducible
objects in $\Mo$ and let the \emph{rank} of $\Mo$ be the cardinal
of $\Irr \Mo$, denoted $\rk \Mo$. If $\Mo$ is a module category of
finite rank then $\Mo$ is a finite direct sum of indecomposable
module categories, since $\rk (\Mo_1\times \Mo_2) = \rk (\Mo_1) +
\rk (\Mo_2)$.

\subsection{Exact module categories over finite tensor categories}\label{subsection:eo}
We are interested in the following class of tensor categories
introduced by Etingof and Ostrik.

\begin{defi}\label{finitetensor}\cite{eo}
Let $\ca$ be a $\ku$-linear category. We shall say that $\ca$ is
\emph{finite} if
\begin{itemize}
\item it has finitely many simple objects; \item each simple
object $X$ has a projective cover $P(X)$; \item the $\Hom$ spaces
are finite-dimensional; \item each object has finite length.
\end{itemize}

A \emph{finite} tensor category is a tensor category $\ca$ such
that the underlying abelian category is finite and the unit object
$\uno$ is  simple.
\end{defi}

\smallbreak  The following definition seems to be part of the
folklore of the subject.

\begin{defi}\label{fusioncat}
A \emph{fusion category} is a finite tensor category $\ca$ such
that the underlying abelian category is semisimple.
\end{defi}

If $H$ is a finite-dimensional Hopf algebra then the category
$\Rep H$ of finite-dimensional representations of $H$ is a finite
tensor category; $\Rep H$ is a fusion category exactly when $H$ is
semisimple.

\smallbreak  It is natural to expect that the study of the module
categories over a finite tensor category $\ca$ would be crucial in
the understanding of $\ca$. As explained in \cite{eo}, one has to
consider a particular class of module categories.

\begin{defi}\label{exactmc} \cite{eo}
A module category $\Mo$ over a finite tensor category $\ca$ is
\emph{exact} if it is finite and for any projective $P\in \ca$ and
any $M\in \Mo$, $P\ot M$ is projective in $\Mo$.
\end{defi}

If $\Mo$ is a semisimple finite module category over a finite
tensor category $\ca$ then it is exact, and the converse is true
if $\ca$ is fusion (tensoring by the unit object).

\begin{rmk}\label{exactmc-frob} \cite[Rmk. 3.7]{eo}
Any exact module category $\Mo$ is a Frobenius category, that is
any projective object of $\Mo$ is injective and vice versa.
\end{rmk}

\begin{rmk}\label{prop:directproduct-exactmc} A direct sum of
finite module categories is exact (resp. semisimple) if and only
if each summand is exact (resp. semisimple). Therefore,  any exact
(resp. finite semisimple) module category over $\ca$ is a finite
direct product of exact (resp. finite semisimple) indecomposable
module categories.
\end{rmk}

\smallbreak

Let $\ca$ be an arbitrary tensor category. A natural way to
produce module categories over $\ca$ is as follows. Let $A$ be an
algebra in $\ca$; then the category $\ca_A$ of right modules in
$\ca$ is a module category over $\ca$.  The category ${}_A\ca_A$
is monoidal with tensor product $\ot_A$ and ${}_A\ca_B$ is a
module category over ${}_A\ca_A$.

\begin{rmk}\label{generalitiesmodcat}
Let $F = (F,b, u):\ca \to \ca'$ be a tensor functor and let $A$ be an algebra
in $\ca$.

(i). $A' := (A', \mu_{A'}, 1_{A'})$, where $A' = F(A)$,
$\mu_{A'} = F(m) b^{-1}$, $1_{A'} = F(1_A) u^{-1}$, is an algebra in $\ca'$.

(ii). If $M\in \ca_A$ is a right $A$-module then $F(M)\in \ca_{A'}$
is a right $A'$-module, with action $\action' = F(\action)b^{-1}$.

(iii). If $\Mo$ is a module category over $\ca'$ with
associativity $m$ and unit $\ell$, then it can be regarded as a
module category over $\ca$ with tensor action $X \ot M := F(X)\ot
M$, $X\in \ca$, $M\in \Mo$, associativity $\widetilde m_{X, Y, M}
:= m_{F(X), F(Y), M}(b^{-1}\ot \id)$ and unit $\ell'_M =
\ell_M(u\ot \id)$.

(iv). Let $\Fc: \ca_A\to \ca'_{A'}$ be the restriction of $F$,
that makes sense by (i) and (ii). Then $(\Fc, c): \ca_A\to
\ca'_{A'}$ is a module functor with respect to the coaction in
(iii), where $c_{X, M} = b_{X,M}$.

(v). If $F$ is an equivalence of tensor categories then
$\Fc$ is an equivalence of module categories.
\end{rmk}

The concept of internal Hom allows to state a converse
of this construction of module categories.

\begin{defi}\label{internalhom}
Let $\Mo$ be a finite module category over $\ca$. Let $M_1, M_2
\in \Mo$. Then the functor $X\mapsto \Hom_{\Mo}(X\ot M_1, M_2)$ is
representable and an object $\uhom(M_1, M_2)$ representing this
functor is called the \emph{internal Hom} of $M_1$ and $M_2$. See
\cite{eo, O1} for details. Thus
$$
\Hom_{\Mo}(X\ot M_1, M_2) \simeq \Hom_{\ca}(X, \uhom(M_1, M_2))
$$
for any $X\in \ca$, $M_1, M_2 \in \Mo$.
\end{defi}

\smallbreak   The ``internal End" $\uend(M) = \uhom(M,M)$ of an
object $M\in \Mo$ is an algebra in $\ca$. The multiplication is
constructed as follows. Denote by $ev_M:\uend(M)\ot M\to M$ the
\emph{evaluation map} \label{evaluation} obtained as the image of
the identity under the isomorphism
$$\Hom_{\ca}(\uend(M),\uend(M))\simeq  \Hom_{\Mo}(\uend(M)\ot M,M).$$
Thus the product $\mu:\uend(M)\ot \uend(M)\to \uend(M)$ is defined
as the image of the map
$$ev_{M} (\id\ot\, ev_{M})\, m_{\uend(M),\uend(M),M}$$ under the
isomorphism
\begin{align*}\Hom_{\ca}(\uend(M)\ot\uend(M),
\uend(M)) \simeq\Hom_{\Mo}((\uend(M)\ot\uend(M)) \ot M, M).
\end{align*}

Recall, on the other hand, that $M\in \Mo$ generates $\Mo$ if for
any $N\in \Mo$, there exists $X\in \ca$ such that $\Hom(X \ot M,
N)\neq 0$. It is known that $W$ generates $\Mo$ iff its simple
subquotients represent all equivalence classes of simple objects
in $\Mo$, where ``equivalence" has the meaning in \cite[Lemma
3.8]{eo}. Hence all simple objects, and \emph{a fortiori} all
non-zero objects, of an indecomposable (finite) module category
are generators.

\begin{teo}\label{exactindecompmc} \cite[Th. 3.17]{eo};
\cite[Th. 3.1]{O2}. Let $\Mo$ be an exact module category over
$\ca$ and let $M\in \Mo$ that generates $\Mo$. Then the functor
$\uhom(M, \underline{\quad}): \Mo \to \ca_{\uend(M)}$ is an
equivalence of module categories.

In particular, if $\ca$ is a finite tensor category and $\Mo$ is
an indecomposable exact module category over $\ca$ then any
non-zero $M\in \Mo$ provides an equivalence of module categories
$\uhom(M, \underline{\quad}): \Mo \to \ca_{\uend(M)}$.
\end{teo}

\begin{rmk}\label{irredimpliesHsimple} \cite[Lemma 4.2]{eo}.
Keep the notation of the Theorem above. The functor $\uhom$ sends
subobjects of $M$ to subobjects of $\ca_{\uend(M)}$; these are the
right ideals of $\uend(M)$. Thus, if $M$ is simple then $\uend(M)$
has no non-zero right ideals.
\end{rmk}

\begin{exa} \cite{Z}, \cite[Prop. 2]{O1}, \cite{D}. Let $G$ be a finite group, $H = \ku G$.
Then all $H$-simple semisimple $H$-module algebras (up to Morita
equivalence) are of the form $\ku G \ot_{H} \End(V)$, where
$H\subseteq G$ is a subgroup and $V$ is a projective
representation of $H$. The left $H$-action on $\End(V)$ is:
$(h.T)(v):=h T(h^{-1} v)$ for all $T\in \End(V), h\in H, v\in V$.
\end{exa}

\subsection{Module categories over finite-dimensional Hopf
algebras} \label{subsection:modcats} Let $H$ be a
finite-dimensional Hopf algebra. Let $K$ be a left $H$-comodule
algebra. Then ${}_K\Mo$ is a left module category over $\Rep H$
via the coaction $\lambda: K \to H\otimes K$. That is, $\ot:\Rep
H\times{}_K\Mo\to {}_K\Mo$ is given by
$$X\ot V:=X\ot_{\ku} V,$$ for $X\in\Rep H$ and $V\in{}_K\Mo$
with  action $k\cdot (x\ot v)=k\_{-1}\cdot x\,\ot \, k\_0\cdot v$,
for all $k\in K$, $x\in X$, $v\in V$.

\smallbreak We say that $K$ is \emph{exact} if $_K\Mo$ is exact,
see Definition \ref{exactmc}. If $K$ is semisimple then it is
exact but the converse is not true, see for example \cite[Th.
4.10]{eo}. However, if $H$ is semisimple then $K$ exact implies
$_K\Mo$ semisimple, as said; hence $K$ semisimple.

\begin{rmk}\label{exactmc-k-frob}
Assume that ${}_K\Mo$ is an exact module category. Then ${}_K\Mo$
is a Frobenius category, that is any projective object of
${}_K\Mo$ is injective and vice versa, see Remark
\ref{exactmc-frob}. In particular $K$ is an injective $K$-module,
that is $K$ is quasi-Frobenius.
\end{rmk}

We next turn to decomposability of ${}_K\Mo$. Recall the
definition of $H$-indecomposable, Definition \ref{hsimple}; it is
clear what ``$H$-decomposable" then means.

\begin{prop}\label{decomposable} The following are equivalent:
\begin{enumerate}
\item The module category ${}_K\Mo$ is  decomposable.

\item The comodule algebra $K$ is $H$-decomposable.
\end{enumerate}
\end{prop}

\begin{proof} (2) $\implies$ (1). Since $I$ and $J$ are $H$-ideals,
the quotient algebras $K/I$ and $K/J$ are left $H$-comodule
algebras. The module categories ${}_{K/I}\Mo$ and ${}_{K/J}\Mo$
are non-trivial submodule categories of ${}_{K}\Mo$. Given $M\in
{}_{K}\Mo$, let $M_1 = \{m\in M: I.m = 0\}$, $M_2 = \{m\in M: J.m
= 0\}$; clearly, $M_1 \in {}_{K/I}\Mo$ and $M_2 \in {}_{K/J}\Mo$.
Decompose $1 = i+ j$, $i\in I$, $j\in J$. Let $m\in M$; then $m =
im + jm$. Since $IJ \subset I\cap J = 0$, $jm \in M_1$, and
similarly $im \in M_2$, thus $M=M_1+M_2$. Also, if $m\in M_1\cap
M_2$, $m=0$. This shows that $M=M_1\oplus M_2$, hence ${}_{K}\Mo
\simeq {}_{K/I}\Mo\times {}_{K/J}\Mo$.

(1) $\implies$ (2). Assume that ${}_{K}\Mo \simeq \Mo_1\times
\Mo_2$, where $\Mo_1$ are $H$-module subcategories of ${}_{K}\Mo$.
If $M\in{}_{K}\Mo$ then there exist $M_1\in\Mo_1$, $M_2\in\Mo_2$
such that $M=M_1\oplus M_2$. If $\varphi: M \to N$ is a morphism
in ${}_{K}\Mo$ then $\varphi(M_1) \subseteq N_1$, $\varphi(M_2)
\subseteq N_2$.

Considering $K\in {}_{K}\Mo$, there exist $J\in \Mo_1$,
$I\in\Mo_2$ such that $K = J\oplus I$. Clearly $I$ and $J$ are
left $H$-ideals of $K$. Let $j\in J$ and let $\eta_j: K\to K$ be
the expansion of $j$, that is $\eta_j(x) = xj$, $x\in K$. Since
$\eta_j$ is a morphism of $K$-modules, $\eta_j(I) \subseteq I$.
Thus, $IJ \subseteq I$ and $I$ is a two-sided ideal. Similarly,
$J$ is a two-sided ideal. Now given an arbitrary $M\in{}_{K}\Mo$
with $M_1\in\Mo_1$, $M_2\in\Mo_2$ such that $M=M_1\oplus M_2$,
then $IM_1 \subset M_1$ and $JM_2 \subset M_2$, by the same
argument applied to the expansion of $m$, $m\in M$. If $I=0$ then
$1\in J$, hence $M= M_2$ and $\Mo_1$ is trivial. Note that these
constructions are inverse to each other.
\end{proof}

\begin{prop}\label{todas}
If $\Mo$ is an indecomposable exact module category over the
finite tensor category $\ca = \Rep H$, then there exists an
indecomposable exact $H$-comodule algebra $K$ such that $\Mo$ is
equivalent to $_K\Mo$ as module categories.
\end{prop}

\pf By Theorem \ref{exactindecompmc} there exists an $H$-module
algebra $R$ such that $\Mo$ is equivalent to $(\Rep H)_R =: \,
_H\Mo_R$. Note that $R^{\op}$ is a $H^{\cop}$-module algebra. Let
$K =R^{\op}\#H^{\cop}$; this is a $H$-comodule algebra.
Explicitly, the multiplication and coaction are respectively given
by
$$
(r\#h)(r'\#t) = (h\_{2}\cdot r') r\# h\_1t, \qquad \lambda(r\# h)
= h\_1 \ot r\# h\_2,
$$
for $r,r'\in R$, $h,t\in H$. Then \begin{align*}
\lambda\big((r\#h)(s\#t)\big) &= h\_1t\_1 \ot (h\_{3}\cdot s) r\#
h\_2t\_1 \\ &= \big(h\_1 \ot r\# h\_2\big)\big(t\_1 \ot s\#
t\_2\big) =\lambda(r\#h)\lambda(s\#t).\end{align*}

Let $\Fc: {}_H\Mo_R\to \,_K\Mo$ be the functor given by $\Fc(M) =
M$, with action $(r\#h)\cdot m = (h\cdot m)\cdot r$, $h\in H$,
$m\in M$, $r\in R$. This is well-defined because the action $M\ot
R \to M$ is a morphism of $H$-modules. Clearly, $\Fc$ is an
equivalence of abelian categories. We claim that $(\Fc, c)$ is an
equivalence of module categories where $c_{X,M}: \Fc(X\otimes
M)\to X\otimes \Fc(M)$ is the identity. Indeed, the only point
that requires some checking is that $c_{X,M}$ is a morphism of
$K$-modules. So, let $X$ be an $H$-module, $M \in {}_H\Mo_R$,
$x\in X$, $m\in M$, $h\in H$, $r\in R$. Then
\begin{align*}
c_{X,M}\big((r\#h)\cdot (x\ot m)\big) &= \big(h\cdot (x\,\ot
m)\big)\cdot r =
h\_1\cdot x\,\ot (h\_2\cdot m)\cdot r; \\
(r\#h)\cdot c_{X,M}(x\,\ot m) &= h\_1\cdot x\,\ot (h\_2\# r)\cdot
m = h\_1\cdot x\,\ot (h\_2\cdot m)\cdot r.
\end{align*}
\epf

\begin{prop}\label{Hsimplealgexact} (i).
If $K$ is an $H$-simple comodule algebra, then $K$ is exact.

(ii). If $K\subseteq H$ is a left coideal subalgebra then $K$ is
exact.
\end{prop}

\pf (i). Let $X$ be a finite-dimensional projective $H$-module and
$M\in {}_K\Mo$. We want to show that $X \otimes M$ is projective;
it is enough to assume that $X = H$. But $H \otimes M \in
{}_K^H\Mo$, hence it is projective as a $K$- module by Thm.
\ref{Skryabin1}. (ii) follows from Prop. \ref{coideal-simple} and
(i). \epf

\subsection{Equivalence of module categories}
Let $H$ be a finite-dimensional Hopf algebra. Let $R$ and $S$ be
left $H$-comodule algebras. We now study module functors between
the module categories ${}_R\mo$ and ${}_S\mo$. For this, we first
recall the following well-known theorem.

\begin{teo}\label{watts}  Let $A$ and $B$ be finite dimensional
algebras. If $F: \,{}_{A}\Mo\to {}_{B}\Mo$ is a right exact
additive functor, then there exists a bimodule $C\in {}_{B}\Mo_A$
and a natural isomorphism $F\simeq C\ot_A \, $. \qed
\end{teo}
\pf The proof goes entirely similar to the proof of \cite[Thm
1]{Wa}; or else it can be deduced from \cite[Thm 2]{Wa}. \epf

\smallbreak To adapt this result to module categories, we
introduce the following notion. First, if $P$ is a
$(S,R)$-bimodule then $H\ot P$ is a $(S,R)$-bimodule by
$$s\cdot(h\ot p)\cdot r = s\_{-1}hr\_{-1} \ot s\_{0}.p.r\_{0},$$ $r\in
R$, $h\in H$, $p\in P$, $s\in S$.

\begin{defi}\label{equivbimodule} An \emph{equivariant $(S,R)$-bimodule}
is a $(S,R)$-bimodule $P$ provided with a left coaction
$\lambda:P\to H\ot_{\ku} P$ that is a morphism of
$(S,R)$-bimodules. Morphisms of equivariant bimodules are defined
in the obvious way. The category of equivariant bimodules is
denoted $\emod$.
\end{defi}

We next prove that the category of module functors $(\Fc,
c):{}_R\mo\to {}_S\mo$ is equivalent to the category of
equivariant $(S,R)$-bimodules.

\smallbreak
\begin{prop}\label{modulefunctors-comalg} There is an equivalence
of categories $\emod \simeq \Hom_{\Rep H}({}_{R}\Mo, {}_{S}\Mo)$.
\end{prop}

\pf Let $P$ be an equivariant $(S,R)$-bimodule. Let $\Fc_P =
\Fc:{}_R\mo\to {}_S\mo$ be the functor defined by $\Fc(V)=P\ot_R
V$. Given $X\in\Rep H$, $V\in {}_R\mo$ we set $c_{X,V}:P\ot_R
(X\ot_{\ku} V)\to X\ot_{\ku} (P\ot_R V)$ by
\begin{equation}\label{comodequi-1} c_{X,V}(p\ot_R x\,\ot
v)=p\_{-1}\cdot x \,\ot\, p\_0\ot_R v, \qquad p\in P, x\in X, v\in
V.
\end{equation}
This is well defined: if $p\in P$, $x\in X$, $v\in V$ and $r\in R$
then
\begin{align*} c_{X,V}(p\cdot r\ot_R x\,\ot v)
&=(p\cdot r)\_{-1}\cdot x \,\ot\, (p\cdot r)\_0\ot_R v
=(p\_{-1}r\_{-1})\cdot x \,\ot\, p\_0\cdot r\_0\ot_R v,\\
c_{X,V}(p\ot_R r\cdot(x\,\ot v)) & = c_{X,V}(p\ot_R r\_{-1}\cdot
x\,\ot r\_{0}\cdot v)= p\_{-1}\cdot (r\_{-1}\cdot x)\,\ot\,
p\_0\ot_R r\_{0}\cdot v.
\end{align*}
Similarly, $c_{X,V}$ is a morphism in ${}_S\mo$. It is an
isomorphism, with inverse given by $c^{-1}_{X,V}(x\,\ot p\ot_R
v)=p\_{0} \ot_R \Ss^{-1}(p\_{-1})\cdot x\,\ot v$, for all $p\in
P$, $x\in X$, $v\in V$; and it is clearly natural. The identities
\eqref{modfunctor1} and \eqref{modfunctor2} are immediate. Hence
$(\Fc, c)$ is a module functor. Furthermore, if $f:P\to Q$ is a
morphism of equivariant $(S,R)$-bimodules then define $\Fc_f:
\Fc_P\to \Fc_Q$ by $\Fc_{f, V} = f\ot \id_V: P\otimes_R V \to
Q\ot_R V$. Equation \eqref{modfunctor3} holds since $f$ is
morphism of $H$-comodules. Thus we have an additive functor $\emod
\to \Hom_{\Rep H}({}_{R}\Mo, {}_{S}\Mo)$.

\smallbreak  Conversely, let $(\Fc, c):{}_R\mo \to {}_S\mo$ be a
module functor. The functor $\Fc$ is exact, \cite[Lemma 3.21]{eo},
hence by Theorem \ref{watts} there exists $P\in {}_S\mo_R$ such
that $\Fc(V)=P\ot_R V$ for all $V\in {}_R\mo$. We define
$\lambda:P\to H\ot_{\ku} P$ by
$$\lambda(p)=c_{H,R}(p\ot 1\ot 1) := p\_{-1}\ot p\_{0}, \qquad p\in P. $$

We first show that $\lambda$ determines $c$. Given $X\in\Rep H$,
$V\in {}_R\mo$, we consider the expansions of $x\in X$, and $v\in
V$, namely $\eta^X_x:H\to X$, $\eta^V_v:R\to V$ given by
$\eta^X_x(h)=h\cdot x$, $\eta^V_v(r)=r\cdot v$ for all $h\in H$,
$r\in R$. The naturality of $c$ implies that the following diagram
is commutative:
\begin{align*}
\begin{CD}
P\ot_R (H\ot_{\ku} R) @>c_{H,R} >> H\ot_{\ku}(P\ot_R R) \\
@V\id\ot\eta^X_x \ot \eta^{V}_v VV    @VV\eta^X_x\,\ot\id\ot\eta^{V}_v V  \\
P\ot_R(X\ot_{\ku} V) @>>c_{X,V}> X\ot_{\ku} (P\ot_R V).
\end{CD}
\end{align*}
Hence for all $p\in P$
\begin{equation}\label{naturalid-c} c_{X,V}(p\ot_R x\,\ot
v)=(\eta^X_x \ot\id_P\ot \eta^{V}_v)\, c_{H,R}(p\ot_R 1\ot 1) =
p\_{-1}\cdot x\,\ot p\_{0}\ot_R v.
\end{equation}

We claim that $P$ is an equivariant bimodule. For this, we first
check that $(P,\lambda)$ is a left $H$-comodule. By naturality of
$c$ in the first variable, the following diagram is commutative:
\begin{align*}
\begin{CD}
P\ot_R (H\ot_{\ku} R) @>c_{H,R} >> H\ot_{\ku}(P\ot_R R) \\
@V\id_P\ot\varepsilon \ot\id VV    @VV\varepsilon\ot\id\ot\id V   \\
P\ot_R(\ku\ot_{\ku} R) @>>c_{\ku,R}> \ku\ot_{\ku} (P\ot_R R).
\end{CD} \end{align*}
The axiom \eqref{modfunctor2} says that $c_{\ku,R}=\id_P$. This
shows that $\lambda$ is counitary.

Again by naturality of $c$ in the first variable, the following
diagram is commutative:
\begin{align}\label{diagram1}
\begin{CD}
P\ot_R (H\ot_{\ku} R) @>c_{H,R} >> H\ot_{\ku}(P\ot_R R) \\
@V\id_P\ot\Delta \ot\id VV    @VV\Delta\ot\id\ot\id V   \\
P\ot_R(H\ot_{\ku} H\ot_{\ku} R) @>>c_{H\ot H,R}> (H\ot_{\ku}
H)\ot_{\ku} (P\ot_R R).
\end{CD}
\end{align}
If $p\in P$ then
\begin{align*} (\Delta\ot\id_P)\lambda(p)&=(\Delta\ot\id_P)\, c_{H,R}(p\ot 1\ot
1)\\
&= c_{H\ot H, R}(\id_P\ot \Delta\ot\id_R)(p\ot 1\ot 1)\\
&=(\id_H\ot c_{H,R}) c_{H,H\ot R}(p\ot 1\ot 1 \ot 1)
\\
&=(\id_H\ot c_{H,R})\lambda(p)\ot 1\ot1\\
&=(\id_H\ot\lambda)\lambda(p).
\end{align*}
The first equality by definition of $\lambda$, the second by
commutativity of diagram \eqref{diagram1}, the third by
\eqref{modfunctor1} and the fourth by \eqref{naturalid-c}. That is
$\lambda$ is coassociative.

\smallbreak We finally  check that $\lambda$ is a morphism of
$(S,R)$-modules. If $p\in P$, $s\in S$, $r\in R$, then
\begin{align*} \lambda(s\cdot p)&=c_{H,R}(s\cdot (p \ot 1\ot 1)) = s\cdot c_{H,R} (p \ot 1\ot 1)
= s\_{-1}p\_{-1} \ot s\_{0}\cdot p\_{-1}\\
\lambda(p\cdot r)&=c_{H,R}(p\cdot r\ot_R (1\ot 1))=c_{H,R}(p\ot_R
r\_{-1}\ot r\_0) = p\_{-1}r\_{-1} \ot p\_{0}\ot_R r\_{0}.
\end{align*} Here, in the first line we have used that $c_{H,R}$ is a morphism of
$S$-modules; and the last equality of the second line comes from
\eqref{naturalid-c}.
 \epf

As a consequence of this result we describe the equivalences of
module categories between ${}_R\mo$ and ${}_S\mo$. Recall that a
\emph{Morita context} for $R $ and $S$, is a collection
$(P,Q,f,g)$ where $P\in {}_S\Mo_R$, $Q\in {}_R\Mo_S$, $f:P\ot_R
Q\stackrel{\simeq}{\to} S$ is an isomorphism of $S$-bimodules and
$g:Q\ot_S P\stackrel{\simeq}{\to} R$ is un isomorphism of
$R$-bimodules such that $ f(p\ot q)p' = pg(q\ot p')$, $g(q\ot p)q'
= qf(p\ot q')$ for all $p,p'\in P$, $q,q'\in Q$.

We shall say that a Morita context $(P,Q,f,g)$ is
\emph{equivariant} if $P$ is an equivariant bimodule. We shall see
that $Q$ turns out to be equivariant too. In this case we shall
say that $R$ and $S$ are \emph{equivariantly Morita equivalent}.

Combining Morita theory with the Proposition
\ref{modulefunctors-comalg}, we get:

\begin{prop}\label{equivcomalg} The equivalences of module categories between
${}_R\mo$ and ${}_S\mo$ are in bijective correspondence with
equivariant Morita contexts for $R$ and $S$.
\end{prop}

\pf Let $F:{}_R\mo\to {}_S\mo$ be an equivalence of module
categories. Then $F$ is, in particular, an equivalence of abelian
$\ku$-linear categories and gives rise to a Morita context
$(P,Q,f,g)$ where $F(M) = P\ot_R M$; furthermore $P$ is an
equivariant bimodule by Proposition \ref{modulefunctors-comalg}.
Conversely, let $(P,Q,f,g)$ be an equivariant Morita context.
Recall that $Q \simeq \Hom_R(P, R)\simeq \Hom_S(P, S)$ and $P
\simeq \Hom_R(Q, R)\simeq \Hom_S(Q, S)$. Then $F:{}_R\mo\to
{}_S\mo$, $F(M) = P\ot_R M$, is an equivalence of $\ku$-linear
categories; its inverse is $G:{}_S\mo\to {}_R\mo$, $G(N) = Q\ot_S
N$ and the natural isomorphisms $\alpha: G\circ F\to
\id_{{}_R\Mo}$, $\beta: F\circ G \to\id_{{}_S\Mo}$,  are given by
\begin{align*}
\alpha&: Q\ot_S P\ot_R M \to M, &  \alpha(q\ot p \ot m) &=
q(p)m, &  M&\in {\,}_R\Mo, \,m\in M; \\
\beta&:  P\ot_R Q\ot_S N \to N, & \beta(p\ot q \ot n) &= q(p)n, &
N&\in {\,}_S\Mo, \, n\in N.
\end{align*}
Here in the first line we have used the identification $Q \simeq
\Hom_R(P, R)$ and in the second, $P \simeq \Hom_S(Q, S)$. We next
consider the left $H$-coaction on $Q \simeq \Hom_R(P, R)$
corresponding to the right $H^*$-action given by
\begin{equation}\label{actionq}
(q \leftharpoonup \gamma)(p) = q(p\leftharpoonup
\Ss^{-1}(\gamma\_{2}))\leftharpoonup\gamma\_{1}, \qquad q\in Q,
p\in P, \gamma\in H^*.
\end{equation}
We claim that $Q$ with this coaction is an equivariant
$(R,S)$-bimodule, which amounts to
\begin{equation}\label{actionq-equivariant}
(rqs) \leftharpoonup \gamma =
(r\leftharpoonup\gamma\_{1})(q\leftharpoonup\gamma\_{2})(s\leftharpoonup\gamma\_{3}),
\qquad q\in Q, r\in R,s\in S, \gamma\in H^*.
\end{equation}
Evaluating both sides at $p\in P$, we have
\begin{align*}
\text{LHS of } \, \eqref{actionq-equivariant} (p) &=
(rqs)(p\leftharpoonup \Ss^{-1}(\gamma\_{2}))\leftharpoonup
\gamma\_{1} = [rq(s(p\leftharpoonup \Ss^{-1}(\gamma\_{2})))]
\leftharpoonup\gamma\_{1}
\\&= (r\leftharpoonup\gamma\_{1})[q(s(p\leftharpoonup
\Ss^{-1}(\gamma\_{3})))\leftharpoonup\gamma\_{2}].
\\ \text{RHS of } \, \eqref{actionq-equivariant} (p) &=
(r\leftharpoonup\gamma\_{1})(q\leftharpoonup\gamma\_{2})((s\leftharpoonup\gamma\_{3})p)
\\ &=
(r\leftharpoonup\gamma\_{1})
[q(((s\leftharpoonup\gamma\_{4})p)\leftharpoonup\Ss^{-1}(\gamma\_{3}))\leftharpoonup\gamma\_{2}]
\\ &=
(r\leftharpoonup\gamma\_{1})[q((s\leftharpoonup\gamma\_{5}\Ss^{-1}(\gamma\_{4}))
(p\leftharpoonup\Ss^{-1}(\gamma\_{3})))\leftharpoonup\gamma\_{2}].
\end{align*}
Thus \eqref{actionq-equivariant} holds. Finally, we show that
$\alpha$ and $\beta$ satisfy \eqref{modfunctor3}. For $\alpha$,
the commutativity of
\begin{gather*}
\xymatrix{ Q\ot_S P\ot_R(X\ot M) \ar[d]_{\alpha_{X\otimes M}}
\ar[rr]^{c_{X,M}}&& X\otimes Q\ot_S P\ot_R M
\ar[d]^{{\id_{X}}\ot\alpha_{M}}\\
X \otimes M \ar[rr]_{\id_{X,M}}&& X\otimes M,} \end{gather*}
\begin{gather*} \xymatrix{ q\ot p\ot x\ot m \ar@{|->}[dd]
\ar@{|->}[r]& q\otimes p\_{-1}\cdot x\otimes  p\_{0} \ot
m\ar@{|->}[r]& q\_{-1}p\_{-1}\cdot x\otimes q\_{0} \ot p\_{0} \ot
m \ar@{|->}[d]
\\
&&q(p)\_{-1}\cdot x \otimes q(p)\_{0}\cdot m \ar@{=}[d]^{\,?}
\\ q(p)(x \otimes m)
\ar@{|->}[rr]&&  q\_{-1}p\_{-1}\cdot x\otimes q\_{0} (p\_{0}) m.}
\end{gather*}
needs the identity $q(p)\_{-1} \otimes q(p)\_{0} =
q\_{-1}p\_{-1}\otimes q\_{0} (p\_{0})$, that follows immediately
from \eqref{actionq}. For $\beta$ the argument is similar once the
agreement of the analog of \eqref{actionq} for the action on $P$
and the original one is shown. \epf

Together with Propositions \ref{decomposable} and \ref{todas},
Proposition \ref{equivcomalg} implies the first approach to the
classification of module categories over $\Rep H$.

\begin{teo}\label{clasif-uno} Indecomposable exact module categories
over a finite dimensional Hopf algebra $H$ are classified by
$H$-indecomposable left comodule algebras up to equivariant Morita
equivalence. \qed
\end{teo}

\begin{lema} (i). Let $P_R$ be a right $R$-module.
Then $\End_R(P_R)$ is a left $H$-comodule algebra via
$\lambda:\End_R(P_R)\to H\ot_{\ku} \End_R(P_R) $, $T\mapsto
T\_{-1}\ot T\_0$, determined by
\begin{equation}\label{h-comod} \langle\alpha, T\_{-1}\rangle\,
T_0(p)=\langle\alpha, T(p\_0)\_{-1}\Ss^{-1}(p\_{-1})\rangle\,
T(p\_0)\_0,\end{equation}  $T\in\End_R(P_R),$ $p\in P$, $\alpha\in
H^*$. Furthermore, $P$ is an equivariant $(S,R)$-bimodule where $S
= \End_R(P_R)$.

(ii). Let $(P,Q,f,g)$ be an equivariant Morita context. The
application $\rho: S\to \End_R(P_R)$, $\rho(s)(p)= s\cdot p$ for
all $s\in S$, $p\in P$, is an isomorphism of $H$-comodules.
\end{lema}

\pf (i). We check that $\lambda $ is well-defined, i. e. that
$T\_{-1}\ot T\_0\in H\ot_{\ku} \End_R(P_R)$. If $\alpha\in H^*$,
$r\in R$, $p\in P$, then
\begin{align*}\langle\alpha, T\_{-1}\rangle\, T_0(r\cdot p)&=
\langle\alpha, T(( p\cdot r)\_0)\_{-1}\Ss^{-1}((p\cdot
r)\_{-1})\rangle\,
T((p\cdot r)\_0)\_0\\
&=\langle\alpha, T(p\_0\cdot r\_0)\_{-1}\Ss^{-1}(
p\_{-1}r\_{-1})\rangle\, T(p\_0\cdot r\_0)\_0\\
&=\langle\alpha,  T( p\_0)\_{-1}r\_{-2}\Ss^{-1}(
p\_{-1}r\_{-1})\rangle\, T( p\_0)\_0\cdot r\_0\\
&=\langle\alpha, T(p\_0)\_{-1}\Ss^{-1}(p\_{-1})\rangle\,
T(p\_0)\_0\cdot r\\
&=\langle\alpha, T\_{-1}\rangle\, T_0(p)\cdot r.
\end{align*}
Here the first equality holds by \eqref{h-comod}, the second
because the coaction $\lambda$ of $P$ is a morphism of
$R$-modules, the third because $T$ is un morphism of right
$R$-modules and the fourth because of properties of $\Ss^{-1}$. It
is immediate that $\lambda$ is coassociative and counitary. It
follows directly from \eqref{h-comod} that $P$ is equivariant.

We now check that $\End_R(P_R)$ is a left $H$-comodule algebra.
Let $\alpha \in H^*$, $T, U\in \End_R(P_R)$ and $p\in P$ then
\begin{align*} &\langle\alpha, T\_{-1} U\_{-1} \rangle\, T\_0
(U\_0(p))=\langle\alpha\_1, T\_{-1}\rangle \langle\alpha\_2,
U\_{-1}\rangle \, T\_0 (U\_0(p))\\
&=\langle\alpha\_1, T\_{-1}\rangle \langle\alpha\_2, U(p\_0)\_{-1}
\Ss^{-1}(p\_{-1})\rangle\, T\_0 (U(p\_0)\_0)\\
&=\langle\alpha,
T(U(p\_0)\_{0})\_{-1}\Ss^{-1}(U(p\_0)\_{-2})U(p\_0)\_{-1}
\Ss^{-1}(p\_{-1})\rangle\, T(U(p\_0)\_{0})\_0\\
&=\langle\alpha, T(U(p\_0))\Ss^{-1}(p\_{-1})\rangle\,
T(U(p\_0)\_0\\
&=\langle\alpha, (TU)\_{-1}\rangle\,  (TU)\_0(p).
\end{align*}
The fourth equality by properties of the antipode. Since this
holds for arbitrary $\alpha\in H^*$ then
$\lambda(TU)=\lambda(T)\lambda(U).$

(ii). By Morita theory, the application $\rho$ above is a linear
isomorphism. Let $s\in S$, $p\in P$ and $\alpha\in H^*$. Then
\begin{align*} \langle \alpha, \rho(s)(p)\_{-1}\rangle\,
\rho(s)(p)\_0&= \langle \alpha, \rho(s)(p\_0)\_{-1}
\Ss^{-1}(p\_{-1})\rangle \, \rho(s)(p\_0)\_0\\
&=\langle \alpha, s\_{-1} p\_{-1} \Ss^{-1}(p\_{-2}) \rangle \,
s\_0\cdot p\_0\\
&=\langle \alpha, s\_{-1} \rangle \, s\_0\cdot p.
\end{align*}
Here the first equality holds by \eqref{h-comod} and the second
because $\lambda$ is a morphism of $S$-modules. Hence
$(\id\ot\rho)\lambda_S= \lambda \rho$. \epf

\begin{rmk}\label{simple=trivialco} The space of coinvariants of
$\End_R(P_R)$ is $\End_R(P_R)^{\co H}=\End^H_R(P_R)$. That is
$T\in \End_R(P_R)$ satisfies $T\_{-1}\ot T\_0=1\ot T$ if and only
if $T$ is an $H$-comodule map. In particular, if $P$ is a simple
object in ${}^H\Mo_R$, then $\End_R(P_R)$ has trivial
coinvariants.

\end{rmk}

\section{Yan-Zhu Stabilizers}\label{section:yz}

In this section, $H$ is a finite-dimensional Hopf algebra.

\subsection{Preliminaries}\label{subsection:modalg}
Let $K$ be a left $H$-comodule algebra, and let $W, U, Z\in
{}_K\Mo$, with corresponding representations $\rho_W: K \to \End
W$, etc. It is convenient to consider the linear map $$\ele = L
\ot \id :H^* \ot \Hom (U, W) \to \Hom(H^* \ot U, H^* \ot W)\simeq
\End (H^*) \ot \Hom(U, W),$$ that is $$\ele(\alpha\ot f)(\beta\ot
u)=\alpha\beta\ot f(u),$$ for all $\alpha, \beta\in H^{*}$,
$f\in\Hom(U, W), u\in U$. We consider the left actions of $H$ on
$H^* \ot \Hom (U, W)$, $H^* \ot U$ and $H^* \ot W$ induced by the
action $\rightharpoonup$ of $H$ on $H^*$ (and trivial on the
second tensorand). In particular, $\Hom (H^* \ot U, H^* \ot W)$
becomes an $H$-module.

\begin{lema}\label{ele}  The map $\ele$ has the following properties:

\begin{itemize}
    \item[(i)] is compatible with
compositions, \emph{i. e.} the following diagram is commutative:

\begin{equation}\label{elecomp}
\begin{CD}
\Hom(H^* \ot U, H^* \ot W) \times \Hom(H^* \ot W, H^* \ot Z)
@>\text{\rm composition}\,>> \Hom(H^* \ot U, H^* \ot Z)
\\ @A\ele\ot  \ele AA @AA \ele A\\
H^* \ot\Hom( U, W) \times H^* \ot\Hom(W, Z) @>\mu \ot\text{\rm
composition}\,>> H^* \ot\Hom(U, Z).
\end{CD}
\end{equation}

\

\smallbreak    \item[(ii)] $\ele:H^* \ot \End W \to \End(H^* \ot
W)$ is a morphism of algebras.

\smallbreak    \item[(iii)] $\ele$ is an injective $H$-module
homomorphism.
\end{itemize}

\end{lema}

\begin{proof}  (i) and (ii) are straightforward.
Clearly $\ele$ is injective; it preserves the action of $H$ since
$L$ does. Indeed, for $h\in H$, $\alpha, \beta\in H^{*}$,
\begin{align*} (h\rightharpoonup L_{\alpha})(\beta)&=h\_1\rightharpoonup
(L_{\alpha}(\Ss(h\_2)\rightharpoonup\beta))=h\_1\rightharpoonup(\alpha(\Ss(h\_2)\rightharpoonup
\beta))\\
&=(h\_1\rightharpoonup\alpha)(h\_2\Ss(h\_3)\rightharpoonup\beta)=L_{h\rightharpoonup
\alpha}(\beta).
\end{align*}
\end{proof}

The discussion above can be carried over for the right regular
representation. Let us consider the map
$$\ere=\id\ot R: \Hom(U,W)\ot H\to \Hom(U\ot H, W\ot H)\cong \Hom(U,W)\ot \End(H)$$
defined by $\ere(f\ot h)(u\ot t)=f(u)\ot th$, where $h,t\in H$,
$f\in\Hom(U, W), u\in U$. We consider the right action of $H^{*}$
on $\Hom(U,W)\ot H$, $U\ot H$, $W\ot H$ induced by the right
action $\leftharpoondown$ of $H^{*}$ on $H$ (and trivial on the
first tensorand), cf. \eqref{accionhomder}. Again $\Hom(U\ot H,
W\ot H)$ is a right $H^{*}$-module.

\begin{lema}\label{ere}  The map $\ere$ has the following properties:
\begin{itemize}
    \item[(i)] $\ere$ is compatible with  compositions,
    \item[(ii)] $\ere:\Hom(U,W)\ot H^{\rm{op}}\to \Hom(U\ot H, W\ot H)$
    is a morphism of algebras, and
    \item[(iii)] $\ere$ is an injective $H^{*}$-module map.
\end{itemize}
\end{lema}

\begin{proof} (i) and (ii) are clear.
The map $\ere$  preserves the action of $H^*$ since $R$ does:
\begin{align*} (R_{h}\leftharpoondown\alpha)(t)&=
(R_h(t\leftharpoondown\Ss^{-1}(\alpha\_2))\leftharpoondown\alpha\_1
=
((t\leftharpoondown\Ss^{-1}(\alpha\_2))h)\leftharpoondown\alpha\_1
\\
&=
(t\leftharpoondown\Ss^{-1}(\alpha\_3)\alpha\_2)(h\leftharpoondown\alpha\_1)
= t(h\leftharpoondown\alpha) = R_{h\leftharpoondown\alpha}(t),
\end{align*}
if $h, t\in H$, $\alpha\in H^{*}$. Here we have used
\eqref{yz-tecnico3}.
\end{proof}

\subsection{Hopf modules}
Recall that a right Hopf module over $H$ is a right $H$-comodule
$M$ with coaction $\rho: M \to M\ot H$ provided with a right
$H$-action such that $\rho$ is a morphism of $H$-modules. Since
$H$ is finite-dimensional, a right Hopf module is the same as a
vector space $M$ provided with a right action of $H$ and a left
action of $H^*$ such that
\begin{align} \label{hopfmod}
\alpha \cdot(m\cdot h) = (\alpha_{(1)} \cdot m)\cdot
(\alpha_{(2)}\rightharpoonup h),
\end{align}
$h\in H$, $\alpha\in H^*$, $m\in M$. If $M$ is a right Hopf module
then the action of $H$ induces an isomorphism $M \simeq M^{\co H}
\ot H$ by the Fundamental Theorem of Hopf modules \cite[Th. 1.9.4,
page 15]{Mo}. Here $M^{\co H} = \{m\in M: \rho(m) = m\ot 1\}=
\{m\in M: \alpha.m = \langle\alpha, 1\rangle m\quad \forall
\alpha\in H^*\}$.

\begin{lema}\label{endhopfmod}
$\End(H^*)$ is a right Hopf module over $H$ with actions
\begin{align} \label{accionHend}
(\alpha \cdot f)(\beta) &= f(\beta \alpha_{(2)})
\Ss^{-1}(\alpha_{(1)}),\\ \label{accionH*end} (f \cdot h)(\beta)
&= f (h\rightharpoonup \beta),
\end{align}
$h\in H$, $\alpha, \beta\in H^*$, $f\in \End(H^*)$. Furthermore
the space of coinvariants $\End(H^*)^{\co H}$ is the image of the
left regular representation $L: H^* \to \End(H^*)$.
\end{lema}

\pf Clearly, \eqref{accionHend} and \eqref{accionH*end} are
respectively a left and a right action; we check \eqref{hopfmod}:
\begin{align*}
\left((\alpha_{(1)} \cdot f)\cdot (\alpha_{(2)}\rightharpoonup
h)\right)(\beta) &= (\alpha_{(1)} \cdot
f)((\alpha_{(2)}\rightharpoonup h) \rightharpoonup \beta)
\\&= f  \left(\left( h_{(1)}
\rightharpoonup \beta\right)\langle \alpha_{(3)}, h_{(2)}\rangle
\alpha_{(2)}\right)\Ss^{-1}(\alpha_{(1)})
\\&= f  \left(\left( h_{(1)}
\rightharpoonup \beta\right)\left( h_{(2)} \rightharpoonup
\alpha_{(2)}\right)\right)\Ss^{-1}(\alpha_{(1)})
\\&= f (h\rightharpoonup (\beta \alpha_{(2)})) \Ss^{-1}(\alpha_{(1)})
\\&= (f\cdot h)(\beta \alpha_{(2)})
\Ss^{-1}(\alpha_{(1)})
\\&= (\alpha \cdot (f\cdot h))(\beta).
\end{align*}

We prove the last statement: given $f\in \End (H^*)$, $f\in \End
(H^*)^{\co H}$ iff $\alpha \cdot f=\alpha (1)f$ for all $\alpha
\in H^*$ iff $f(\beta\alpha) =f(\beta)\alpha$ for all $\alpha,
\beta \in H^*$ iff $f =L_\gamma$ for $\gamma= f(1)\in H^*$. \epf

If $K$ is a finite-dimensional Hopf algebra, a right Hopf module
over $K^{\rm{cop}}$ is the same as a vector space $M$ provided
with right actions of $K$ and  $K^*$ such that
\begin{align} \label{hopfmod-cop}
 (m\cdot k)\cdot \gamma = (m\cdot \gamma_{(1)})\cdot
(k\leftharpoonup \gamma_{(2)}),
\end{align}
$k\in K$, $\gamma\in K^*$, where all notations are in terms of
$K$. We shall also need the following.

\begin{lema}\label{dualhopfmodule} $\End(H)$ is a right Hopf
module over $H^{*\rm{cop}}$ with actions
\begin{align}\label{dualhopfmodule1} (f\cdot h)(t)&=\Ss(h\_1)f(h\_2t),
\\\label{dualhopfmodule2}
(f\cdot\alpha)(t)&=f(\alpha \rightharpoonup t),
\end{align}
$h, t\in H$, $\alpha\in H^{*}$, $f\in \End(H)$. Moreover the space
of coinvariants $\End(H)^{\rm{co}\, H^{*\rm{cop}}}$ is the image
of the right regular representation $R:H\to \End(H)$. \qed
\end{lema}

\subsection{The Heisenberg double}\label{subsection:heis}
Recall that the Heisenberg double $\Hc(H^*)$ of the Hopf algebra
$H^*$ is the vector space $H^* \ot H$ with the multiplication $
(\alpha\ot h)(\alpha' \ot h') = \alpha (h_{(1)}\rightharpoonup
\alpha')\ot h_{(2)}h'$. Here, and in the next proposition, $h,h',
t, u\in H$, $\alpha, \alpha', \beta \in H^*$.

\smallbreak
\begin{prop}\label{heisenberg} (i).
There is an isomorphism of algebras $\Psi_1:\Hc(H^*) \to\End H$
given by $\Psi_1(\alpha \ot h)(t) = \alpha\rightharpoondown (ht)$.

\smallbreak  (ii). There is an isomorphism of algebras
$\Psi_2:\Hc(H^*) \to \End (H^*)$ given by\newline  $\Psi_2(\alpha
\ot h)(\beta) = \alpha (h\rightharpoonup\beta)$.

\smallbreak  (iii). The isomorphism of algebras $\Psi:\End H \to
\End (H^*)$, $\Psi = \Psi_2\Psi_1^{-1}$ satisfies
\begin{align}
\label{uno}
\Psi(\elt_{\alpha}) &= L_{\alpha},\\
\label{dos}\Psi(L(H)) &= \el(H),\\
\label{tress} \Psi(\el(H^*))&=R(H^*),\\
\label{cuatros}\Psi(R(H))&=\elt(H).
\end{align}
\end{prop}

\pf We check (i):
\begin{align*}
\Psi_1(\alpha\ot h)\Psi_1(\alpha'\ot h')(t) &= \Psi_1(\alpha\ot
h)(\alpha'\rightharpoondown (h't)) = \alpha\rightharpoondown
\big(h (\alpha'\rightharpoondown (h't))\big);\\
\Psi_1\big((\alpha\ot h)(\alpha'\ot h')\big)(t) &= \Psi_1(\alpha
(h_{(1)}\rightharpoonup \alpha')\ot h_{(2)}h')(t) = \big(\alpha
(h_{(1)}\rightharpoonup \alpha')\big)\rightharpoondown
(h_{(2)}h't)
\\ &= \alpha\rightharpoondown \big((h_{(1)}\rightharpoonup \alpha')
\rightharpoondown (h_{(2)}h't)\big).
\end{align*}
Then (i) follows from the identity $ h (\alpha'\rightharpoondown
u)= (h_{(1)}\rightharpoonup \alpha') \rightharpoondown
(h_{(2)}u)$, that we prove next:
\begin{align*}
(h_{(1)}\rightharpoonup \alpha') \rightharpoondown (h_{(2)}u) &=
\langle h_{(1)}\rightharpoonup \alpha', \Ss^{-1}(h\_2u\_1)\rangle
\, h\_3u\_2
\\ &= \langle \alpha', \Ss^{-1}(u\_1)\Ss^{-1}(h\_2) h\_1\rangle
\, h\_3u\_2
\\ &= h\langle \alpha', \Ss^{-1}(u\_1)\rangle u\_2 =
h (\alpha'\rightharpoondown u).
\end{align*}
It is well-known that the Heisenberg double is a simple algebra,
hence $\Psi_1$ is an isomorphism by a dimension argument. We check
(ii):
\begin{align*}
\Psi_2(\alpha\ot h)\Psi_2(\alpha'\ot h')(\beta) &=
\Psi_2(\alpha\ot h)(\alpha' (h'\rightharpoonup\beta)) =
\alpha (h\rightharpoonup(\alpha' (h'\rightharpoonup\beta)))\\
&= \alpha (h\_1\rightharpoonup \alpha')
(h\_2h'\rightharpoonup\beta) = \Psi_2(\alpha
(h_{(1)}\rightharpoonup \alpha')\ot h_{(2)}h')(\beta) \\ &=
\Psi_2\big((\alpha\ot h)(\alpha'\ot h')\big)(\beta).
\end{align*}

Since $\Psi_1(\alpha \ot 1) = \elt_{\alpha}$ and $\Psi_2(\alpha
\ot 1) = L_{\alpha}$, \eqref{uno} holds. Since $\Psi_1(\varepsilon
\ot h) = L_h$ and $\Psi_2(\varepsilon \ot h) = \underline L_h$,
\eqref{dos} holds. Equation \eqref{tress} follows from \eqref{uno}
and $\elt(H^*)'=\el(H^*)$. Similarly \eqref{cuatros} follows from
\eqref{dos} and $\el(H)'=\elt(H)$. Here $A'$ means the centralizer
of a subalgebra $A$ of $\End V$, see page \pageref{dualproducto}
below. \epf

\subsection{Definition of Yan-Zhu Stabilizers}\label{subsection:stabilizers}
Let $K$ be a left $H$-comodule algebra. We recall the construction
of the stabilizer from \cite{YZ}. Let us consider $H^*$ as an
$H$-module via $\rightharpoondown$, see \eqref{arponbajo}; and
correspondingly $H^* \ot W$ as a $K$-module via $\lambda$. That
is,
\begin{equation}\label{repre}
k\cdot(\beta \ot w) = k_{(-1)} \rightharpoondown \beta\, \ot
k_{(0)}\cdot w,
\end{equation}
$k\in K$, $\beta\in H^*$, $w\in W$. Recall the map $\ele$
considered in subsection \ref{subsection:modalg}.

\smallbreak
\begin{defi}\cite{YZ} The {\it Yan-Zhu stabilizer} of the $K$-modules
$W$ and $U$ is
$$\St_K(U, W):= \Hom_K(H^* \ot U,H^* \ot W)
\cap \ele\big(H^* \ot \Hom (U, W)\big) \subset \Hom(H^* \ot U,H^*
\ot W). $$ In particular, the {\it Yan-Zhu  stabilizer} of the
$K$-module $W$ is
$$\St_K(W):= \St_K(W, W) = \End_K(H^* \ot W) \cap \ele (H^* \ot \End W). $$
\end{defi}

\smallbreak  The algebra $\St_K(W)$ can be identified with a
subalgebra of $H^*\ot \End(W)$, since $\ele$ is injective.
Similarly, $\St_K(U,W)$ can be identified with a subspace of
$H^*\ot \Hom(U,W)$. We shall do this  without further notice.

\begin{prop}\label{stab=rcomalg} $\St_K(W)$ is a right
$H^*$-comodule algebra and $W$ is a left $\St_K(W)$-module.

\end{prop}

\pf By \eqref{elecomp}, the composition induces a map
\begin{equation}\label{stabcomp}
\St_K(U, W) \times \St_K(W, Z) \longrightarrow\St_K(U, Z).
\end{equation}
In particular, $\St_K(W)$ is a subalgebra of $\End(H^*\ot W)$.
Since $\End(H^*\ot W)$ is an $H$-module algebra and $\ele (H^* \ot
\End W)$ is an $H$-submodule, it remains only to show that
$\End_K(H^* \ot W)$ is also an $H$-submodule. But, more generally,
$\Hom_K(H^* \ot U,H^* \ot W)$ is an $H$-submodule of $\Hom(H^* \ot
U,H^* \ot W)$ because $\rightharpoonup$ and $\rightharpoondown$
commute. Hence $\St_K(W)$ is a left $H$-module algebra, thus a
right a $H^*$-comodule algebra. The vector space $W$ is a left
$\St_K(W)$-module, with representation given by the composition
\newline $\begin{CD} \St_K(W) @>>> H^*\ot \End(W)@>\varepsilon\ot \id>>
\ku\ot \End(W) = \End (W).
\end{CD}$
\epf

 The following Lemma will be useful later.

\begin{lema} Let $\sum_i L_{\alpha_i}\ot f_i \in
\ele\left(H^{*}\ot\Hom(U,W)\right)$; we shall omit the summation
symbol. Then $L_{\alpha_i}\ot f_i\in \Hom_K(H^* \ot U, H^* \ot W)$
if and only if
\begin{equation}\label{stab1} \sum_i \langle\alpha_i, t\rangle \, f_i(k\cdot w)=
\sum_i \langle\alpha_i, \Ss^{-1}(k\_{-1})t\rangle \,\, k\_0\cdot
f_i( w),
\end{equation}
for all $t\in H$, $w\in W$, $k\in K$.
\end{lema}

In other words, \eqref{stab1} characterizes when $\sum_i
L_{\alpha_i}\ot f_i\in \St_K(U, W)$.

\begin{proof}
If $L_{\alpha_i} \ot f_i\in \Hom_K(H^* \ot U, H^* \ot W)$ then
\begin{equation}\label{stab1bis} (L_{\alpha_i}\ot
f_i)(k\cdot (\beta\ot w))= k\cdot (L_{\alpha_i}\ot f_i)(\beta\ot
w)
\end{equation}
for all $\beta\in H^{*}, k\in K, w\in W$. Equation
\eqref{stab1bis} translates into
\begin{equation}\label{stab1ter}\alpha_i(k\_{-1}\rightharpoondown\beta)\ot
f_i(k\_0 \cdot w)=k\_{-1}\rightharpoondown (\alpha_i\beta)\ot
k\_0\cdot f_i(w).
\end{equation}
Evaluating \eqref{stab1ter} in $t\in H$ and choosing
$\beta=\varepsilon$ we get \eqref{stab1}. Conversely,
\eqref{stab1ter} follows from \eqref{stab1} using that
$h\rightharpoondown (\alpha\beta) = (h\_2\rightharpoondown
\alpha)\,(h\_1\rightharpoondown \beta)$, for $h\in H$, $\alpha,
\beta\in H^*$.
\end{proof}

As a consequence of the above characterization of elements in the
stabilizer, we have the following result.

\begin{cor} There is an isomorphism $\St_K(U,W)^{\co H^{*}}\cong \Hom_K(U,W)$. In
particular the stabilizer $\St_K(W)$ has trivial coinvariants if
$W$ is an irreducible $K$-module.

\end{cor}

\begin{proof} The maps
$$\phi:\St_K(U,W)^{\co H^{*}}\to \Hom_K(U,W) \quad \text{and}
\quad\psi:\Hom_K(U,W)\to \St_K(U,W)^{\co H^{*}}$$ given by
$\phi(\alpha_i\ot f_i)=\langle\alpha_i, 1\rangle f_i$, $
\psi(f)=\epsilon\ot f$, are the desired isomorphisms.
\end{proof}

We now state another characterization of the stabilizer
$\St_K(U,W)$ given in \cite{YZ} in terms of Hopf modules. We
consider $\End(H^*)\ot \Hom (U,W)$ as right Hopf module over $H$
with structure concentrated in the first tensorand, \emph{cf.}
Lemma \ref{endhopfmod}. We stress that these are not the same
actions as before.

\begin{prop}\label{stabmodhopf} Keep the notation above.

\begin{enumerate}
\item $\Hom_K(H^* \ot U,H^* \ot W)$ is a Hopf submodule
of $\End(H^*)\ot \Hom (U,W)$,
\item $\Hom_K(H^* \ot U,H^* \ot W)^{\co H} = \St_K(U,W)$ and
$$\Hom_K(H^* \ot U,H^* \ot W) = \St_K(U,W)\circ (\el(H)\ot\id) \simeq \St_K(U,W) \ot H.$$
Here $\circ$ means composition. In particular,
\begin{equation}\label{aplhopfmod}
\End_K(H^* \ot W) = \St_K(W)\circ (\el(H)\ot\id).
\end{equation}
\end{enumerate}
\end{prop}

\pf (1). We have to check that $\Hom_K(H^* \ot U,H^* \ot W)$ is
stable under the actions induced by \eqref{accionHend},
\eqref{accionH*end}. Let $k\in K$, $\alpha, \beta\in H^*$, $u\in
U$, $h\in H$. Let also $\sum_i f_i\ot T_i\in \Hom_K(H^* \ot U,H^*
\ot W)$, where $f_i\in \End(H^*)$, $T_i\in \Hom(W, U)$; for
simplicity we omit the summation symbol in the following. Then
\begin{align*}
\big(\alpha \cdot (f_i\ot T_i)\big)\big(k &\cdot(\beta\ot w)\big) =
(\alpha \cdot f_i)(k_{(-1)}\rightharpoondown \beta) \ot
T_i(k_{(0)}\cdot u)
\\ &=  f_i((k_{(-1)}\rightharpoondown \beta) \alpha_{(2)})
\Ss^{-1}(\alpha_{(1)})\ot T_i(k_{(0)}\cdot u)
\\ &=  f_i((k_{(-1)}\leftharpoonup\alpha_{(2)})
\rightharpoondown \beta \alpha_{(3)})
\Ss^{-1}(\alpha_{(1)})\ot T_i(k_{(0)}\cdot u)
\\ &= f_i(\langle k_{(-2)}, \alpha_{(2)}\rangle k_{(-1)}
\rightharpoondown \beta \alpha_{(3)}) \Ss^{-1}(\alpha_{(1)})\ot
T_i(k_{(0)}\cdot u)
\\ &=  \big( (k_{(-1)}\leftharpoonup\alpha_{(2)})
\rightharpoondown f_i( \beta \alpha_{(3)}) \big)
\Ss^{-1}(\alpha_{(1)})\ot k_{(0)}\cdot T_i(u)
\\ &=  \big( (k_{(-1)}\leftharpoonup\alpha_{(3)})
\leftharpoonup \Ss^{-1}(\alpha_{(2)})
\big)\rightharpoondown \big(f_i( \beta \alpha_{(4)})
\Ss^{-1}(\alpha_{(1)})\big)\ot k_{(0)}\cdot T_i(u)
\\ &=  k_{(-1)}\rightharpoondown \big(f_i( \beta \alpha_{(2)})
\Ss^{-1}(\alpha_{(1)})\big)\ot k_{(0)}\cdot T_i(u)
\\ &= k \cdot\Big( \big(\alpha \cdot (f_i\ot T_i)\big)(\beta\ot
w)\Big).
\end{align*}
Here the first two equalities, the fourth and the last are by
definitions; the third and the sixth by \eqref{yz-tecnico}; the
fifth, because $\sum_i f_i\ot T_i\in \Hom_K(H^* \ot U,H^* \ot W)$;
the seventh, by elementary properties of the antipode. Next,
$(f_i\ot T_i)\cdot h$ preserves the $K$-action since
$\rightharpoondown$ and $\rightharpoonup$ commute. Thus, (1)
holds.

(2). We have
\begin{align*}
\Hom_K(H^* \ot U,H^* \ot W)^{\co H} &=\Hom_K(H^* \ot U,H^* \ot W)
\cap \big(\End(H^*)\ot \Hom (U,W)\big)^{\co H} \\&= \Hom_K(H^* \ot U,H^* \ot W)
\cap \End(H^*)^{\co H}\ot \Hom (U,W)
 \\&= \Hom_K(H^* \ot U,H^* \ot W)\cap L(H^*)\ot \Hom (U,W)
\\&= \St_K(U,W).
\end{align*}
The last statement follows from the Fundamental Theorem of Hopf
modules. \epf

\subsection{Yan-Zhu duality}\label{subsection:yzduality}
Let now $S$ be a \emph{right} $H^*$-comodule algebra. Let $V,X,Y$
be left $S$-modules. We adapt the construction of the Yan-Zhu
stabilizer in the new setting. We consider $H$ as a left
$H^*$-module via $\rightharpoondown$ and $V \ot H$ as $S$-module
via the right coaction, that is
\begin{equation}\label{tres}
s\cdot (v\ot h)=s\_0\cdot v\ot s\_1\rightharpoondown h,
\end{equation}
where $s\in S, v\in V$ and $\delta: S\to S\ot H$,
$\delta(s)=s\_0\ot s\_1$. Recall the map $\ere$ considered in
subsection \ref{subsection:modalg}. Then the {\it Yan-Zhu
stabilizer} of the $S$-modules $V$ and $Y$ is
$$\St_S(V, Y):= \Hom_S( V\ot H,  Y\ot H) \cap \ere (\Hom(V, Y)\ot  H). $$
In particular, the {\it Yan-Zhu  stabilizer} of the $S$-module $V$
is $\St_S(V):= \End_S( V\ot H) \cap \ere (\End (V)\ot H)$. We
consider $ \Hom (V, Y)\ot\End (H)$ as right Hopf module over
$H^{*\rm{cop}}$ with structure concentrated in the second
tensorand, \emph{cf.} Lemma \ref{dualhopfmodule}.  Adapting the
proofs of Propositions \ref{stabmodhopf} and \ref{stab=rcomalg},
we have:

\begin{prop}\label{dual1} \

\begin{enumerate}
\item $\Hom_S(V \ot H, Y \ot H)$ is a Hopf submodule of $\Hom (V,
Y) \ot \End (H)$,

 \item $\Hom_S(V \ot H, Y \ot H)^{\co H^{*\rm{cop}}} =
\St_S(V,Y)$,

 \item $\Hom_S(V \ot H, Y \ot H) = \St_S(V,Y) \circ
(\id\ot \el(H^*)) \simeq \St_S(V,Y) \ot H^*$.
\end{enumerate}
\end{prop}

\pf (1). We have to check that $\Hom_S(V\ot H,Y\ot H)$ is stable
under the actions induced by \eqref{dualhopfmodule1},
\eqref{dualhopfmodule2}. Let $s\in S$, $\alpha\in H^*$, $v\in V$,
$h, t\in H$; and let $\sum_j T_j\ot f_j\in \Hom_S(V \ot H,Y \ot
H)$, where $f_j\in \End(H)$, $T_j\in \Hom(V, Y)$; for simplicity
we omit the summation symbol in the following. First, $(T_j\ot
f_j)\cdot \alpha$ preserves the $S$-action since
$\rightharpoondown$ and $\rightharpoonup$ commute. Next,
\begin{align*}
s \cdot ((T_j\ot f_j)\cdot h )(v\ot t)&= s_{(0)}\cdot T_j(v) \ot
s_{(1)}\rightharpoondown (f_j\cdot h)(t)
\\ &=s_{(0)}\cdot T_j(v) \ot s_{(1)}\rightharpoondown
\big(\Ss(h\_1)f_j( h\_2t)\big)
\\ &=s_{(0)}\cdot T_j(v) \ot \big(s_{(2)}\rightharpoondown \Ss(h\_1)\big)
\big(s_{(1)}\rightharpoondown f_j( h\_2t)\big)
\\ &= T_j(s_{(0)}\cdot v) \ot \big(s_{(2)}\rightharpoondown \Ss(h\_1)\big)
 f_j(s_{(1)}\rightharpoondown (h\_2t))
 \\ &= T_j(s_{(0)}\cdot v) \ot \big(s_{(3)}\rightharpoondown \Ss(h\_1)\big)
 f_j((s_{(2)}\rightharpoondown h\_2)(s_{(1)}\rightharpoondown t))
\\ &=  T_j(s_{(0)}\cdot v) \ot \Ss(h\_1)f_j(h\_2(s_{(1)}\rightharpoondown
t))
\\ &=
T_j(s_{(0)}\cdot v) \ot (f_j\cdot h)(s_{(1)}\rightharpoondown t)
\\ &=\big((T_j\ot f_j)\cdot h\big)\big(s \cdot(v\ot t)\big).
\end{align*}
Here the only equality that needs explanation is the sixth, which
is based on the following:
\begin{align*}
\alpha\_2\rightharpoondown \Ss(h\_1)\ot \alpha\_1\rightharpoondown
h\_2&= \langle \alpha\_2, \Ss^{-1}(\Ss(h\_2))\rangle \Ss(h\_1) \ot
\langle \alpha\_1, \Ss^{-1}(h\_3)\rangle \Ss(h\_4)
\\ &= \langle\alpha, 1\rangle\Ss(h\_1)\ot h\_2.
\end{align*}
Thus, (1) holds. The proof of (2) is similar the the proof of
proposition \ref{stabmodhopf} part (2), and (3) follows again from
(2) and the Fundamental Theorem of Hopf modules. \epf

\smallbreak  To state the next result (Yan-Zhu duality), we use
the following notation: if $A$ is a subspace of $\End(W)$ then $A'
= \Cent_{\End (W)}(A)$ is the centralizer of $A$ in $\End(W)$.
Clearly:
\begin{equation}\label{dualproducto}
\text{ If } A = B\circ C \text{ and } 1\in B \cap C \text{ then }
A' = B'\cap C'.
\end{equation}

If $\phi:\End(W)\to \End(V)$ is an algebra isomorphism, then
$\phi(A')=\phi(A)'$. If $A$ is an algebra and $\rho:A \to \End(
W)$ is a representation then $\rho(A)'$ is nothing but
$\End_A(W)$.

Let us fix  a right $H^*$-comodule algebra $S$ and a left
$S$-module $W$; therefore $W$ is a left $\St_S(W)$-module by
proposition \ref{stab=rcomalg}.

\begin{prop}\label{prop:yzduality} There is an isomorphism
of left $H$-module algebras $$ \St_{\St_S(W)}(W) \simeq
\big(\rho_{ W\ot H}(S)\big)'',$$ where $\rho_{W\ot H}$ is the
representation of $S$ explained in \eqref{tres}.
\end{prop}

\pf Since $\St_S(W)$ is a left $H$-comodule algebra and $W$ is a
$\St_S(W)$-module-- see proposition \ref{stab=rcomalg}-- there is
a representation $\rho_{H^*\ot W}:\St_S(W) \to \End (H^*\ot W)$,
given by \eqref{repre}. Recall the isomorphism of algebras $\Psi:
\End H \to \End (H^*)$ given in proposition \ref{heisenberg}. We
claim that
\begin{align}\label{cuatro}(\Psi^{-1}\ot \id)(\rho_{H^*\ot W}\St_S(W)) =
\St_S(W).\end{align} This follows from the definitions and
\eqref{cuatros}. Let now $\Upsilon = (\id\ot\Psi^{-1} )\tau: \End
(H^*\ot W) \to \End (W\ot H)$, where $\tau: \End (H^*\ot W) \to
\End (W\ot H^*)$ is the usual transposition. Then
\begin{align*}
\big(\rho_{W\ot H}(S)\big)'' &= (\End_S( W\ot H))'
\\ &= \big(\St_S(W)\circ ( \id\ot \el(H^*))\big)'\\
&=\St_S(W)'\cap ( \id\ot \el(H^*))\big)'\\
&=\Upsilon\Upsilon^{-1}\left(\St_S(W)'\right)\cap
\Upsilon\Upsilon^{-1}\left(( \id\ot \el(H^*)\big)''\right)\\
&=\Upsilon\left(\Cent_{\End(H^*\ot W)}
\Upsilon^{-1}\St_S(W)\right)
\cap \Upsilon\left(\Upsilon^{-1}(\id\ot \el(H^*))\right)'\\
&=\Upsilon\left(\End_{\Upsilon^{-1}\St_S(W)} (H^*\ot W)\right)\cap
\Upsilon\left(\id\ot R(H^*))'\right)\\
&=\Upsilon\left(\End_{\St_S(W)}(H^*\ot W)\cap
\ele(H^*\End(W))\right)\\
&=\Upsilon\left(\St_{\St_S(W)}(W)\right).
\end{align*}
Here the second equality holds by Proposition \ref{dual1} (3); the
third by \eqref{dualproducto}; the fourth, because $\Upsilon$ is
an algebra isomorphism; the fifth is a restatement; the sixth is
by \eqref{tress};  the seventh follows from \eqref{cuatro} and the
last is by definition. The left $H$-action on $\End(W\ot H)$ is
given by $h\cdot (T\ot f)= T\ot L_{h\_1} f L_{\Ss(h\_2)}$, for all
$h\in H$, $T\in \End(W)$, $f\in \End H$. The map $\Upsilon $ is an
$H$-module map since for all $h\in H$, $\alpha\in H^*$
\begin{align*} \Psi^{-1}(h\cdot L_{\alpha})= \Psi^{-1}(L_{h\rightharpoonup\alpha})
=\elt_{h\rightharpoonup\alpha}=L_{h\_1}\elt_{\alpha}L_{\Ss(h\_2)}=h\cdot
\Psi^{-1}(L_{\alpha}).
\end{align*}
The second equality by \eqref{uno}, the third follows from the
identity $$(h\rightharpoonup \alpha)\rightharpoondown t=
h\_1(\alpha\rightharpoondown(\Ss(h\_2))t), \quad h,t \in H,
\alpha\in H^*.$$ \epf

If $A$ is a quasi-Frobenius algebra and $M$ is a faithful finitely
generated $A$-module then $(A;M)$ has the double centralizer
property, see \cite[Th. 15.6]{CR}. In view of this, and as a
consequence of proposition \ref{prop:yzduality}, we have:

\begin{cor}\label{cor:yzduality} Assume that
\begin{enumerate}
    \item $S$ is a quasi-Frobenius algebra, and
    \item $\rho_{ W\ot H}: S\to \End(W\ot H)$ is injective.
\end{enumerate}
Then $\St_{\St_S(W)}(W)$ is isomorphic to $S$ as $H^*$-comodule
algebras. \qed
\end{cor}

Here is the main result of this Section, proved in \cite{YZ}
assuming that $H$ and $K$ are semisimple.

\begin{teo}\label{cor:yzduality-simple} Let $S$ be an $H$-simple left module algebra.
Then $\St_{\St_S(W)}(W)$ is isomorphic to $S$ as $H$-module
algebras. \qed
\end{teo}

\pf We need to analyze the hypotheses in corollary
\ref{cor:yzduality}. The injectivity of $\rho_{ W\ot H}$ is
disposed with Lemma \ref{inyectiva}. Now $S$ is exact by
Proposition \ref{Hsimplealgexact}, which follows from Skryabin's
Theorem \ref{Skryabin1}. Hence, $S$ is quasi-Frobenius by Remark
\ref{exactmc-frob}. \epf

\subsection{Dimension of Yan-Zhu stabilizers}
We prove in this subsection a formula on the dimension of the
stabilizers, generalizing \cite[Cor. 2.8]{Z}. We begin by a
technical Lemma. Let $K$ be a finite-dimensional $H$-simple left
$H$-comodule algebra. Recall the action of $H$ on $H^*$ given by
$\rightharpoondown$.

\begin{lema}\label{tensor-dual} If $W$ is a left $K$-module then
$H^*\ot W$ is an object in ${}_K^H\Mo.$
\end{lema}
\pf The left $H$-comodule structure $\delta: H^*\ot W\to H\ot
H^*\ot W$ is given as follows. If $\alpha \in H^*$, $w\in W$ then
$\delta(\alpha\ot w)=\alpha\_{-1}\ot\, \alpha\_0\ot\, w$ where
$\langle \beta, \alpha\_{-1}\rangle\, \alpha\_0 =\alpha\,
\Ss^2(\beta)$, for all $\beta\in H^*$. This left $H$-coaction
corresponds to the right $H^*$-action given by the multiplication
composed with $\Ss^2$. Let us verify that the map $\delta$ is a
$K$-module map. Let $k\in K$, $\alpha, \beta \in H^*$ and $w\in
W$. Then
\begin{equation}\label{comod11} \delta(k\cdot ( \alpha\ot w))=(k\_{-1}\rightharpoondown
\alpha)\_{-1}\ot \, (k\_{-1}\rightharpoondown \alpha)\_0\ot\,
k\_0\cdot w \end{equation}  Evaluating  $\beta$ on the first
tensorand of \eqref{comod11} we obtain
\begin{align}\label{com-d1} \langle \beta,(k\_{-1}\rightharpoondown
\alpha)\_{-1}\rangle\,(k\_{-1}\rightharpoondown \alpha)\_0\ot\,
k\_0\cdot w=(k\_{-1}\rightharpoondown \alpha)\, \Ss^2(\beta)\ot\,
k\_0\cdot w.
\end{align}
On the other hand
\begin{equation}\label{comod12} k\cdot \delta( \alpha\ot w)= k\_{-2}\alpha\_{-1}\ot\,
k\_{-1}\rightharpoondown \alpha\_0\ot\, k\_0\cdot w.\end{equation}
Again, evaluating $\beta$ on the first tensorand of
\eqref{comod12} we obtain
\begin{align*} &\langle \beta,
k\_{-2}\alpha\_{-1}\rangle\, k\_{-1}\rightharpoondown
\alpha\_0\ot\, k\_0\cdot w\\& =\langle \beta\_1 , k\_{-2}\rangle
\langle \beta\_2,\alpha\_{-1}\rangle k\_{-1}\rightharpoondown
\alpha\_0\ot\, k\_0\cdot w \\
&=\langle \beta\_1 , k\_{-2}\rangle \left(k\_{-1}\rightharpoondown
(\alpha\, \Ss^2(\beta\_2)\right)\ot\, k\_0\cdot w \\
&=\left( k\_{-1}\leftharpoonup
\Ss^{-2}(\Ss^2(\beta)\_1)\right)\rightharpoondown (\alpha\,
\Ss^2(\beta)\_2)\ot\, k\_0\cdot w \\
&= (k\_{-1}\rightharpoondown \alpha)\Ss^2(\beta)\ot\, k\_0\cdot w.
\end{align*}
The last equality follows from \eqref{yz-tecnico}. Since $\beta$
is arbitrary $\delta(k\cdot ( \alpha\ot w))=k\cdot \delta(
\alpha\ot w)$. \epf

The next formula was obtained in \cite[2.8]{Z} for $H$ a
semisimple Hopf algebra and $U=W$.

\begin{prop}\label{dim-stab} Let $K$ be an $H$-simple left $H$-comodule
algebra and $U$, $W$ two left $K$-modules. Then
\begin{equation}\label{dimension}\dim K \dim \St_K(U,W)=\dim U \dim W \dim H.\end{equation}
\end{prop}

\pf Let $M=H^*\ot W$, $N=H^*\ot U$. By Theorem \ref{freeness} and
Lemma \ref{tensor-dual}, there exists $t, s\in \Na$ such that
$M^t$ and $N^s$ are $K$-free, say $M^t \simeq K^d$ and $N^s\simeq
K^c$ as left $K$-modules for some natural numbers $d, c$. Hence
\begin{align}\label{dim11} t \dim H \dim W &= d \dim K,
\\\label{dim111} s \dim H \dim U &= c \dim K.
\end{align}
Proposition \ref{stabmodhopf} (2) implies that
\begin{equation}\label{dim12} \dim \St_K(U,W) \dim H= \dim
\Hom_K(N,M).
\end{equation}
Since $M^t \simeq K^d$, $N^s\simeq K^c$, it follows that
$\Hom_K(N^s,M^t)\simeq \Hom_K(K^c,K^d)$. Therefore $t s
\dim\Hom_K(N,M)\, =d c\dim K$. This equation combined with
\eqref{dim12} implies that $\displaystyle\dim \St_K(W) \dim H=
\frac{dc\dim K }{ts}$. The result now follows from \eqref{dim11}
and \eqref{dim111}. \epf

As an immediate consequence of \eqref{dimension} we obtain the
following variation of \cite[Prop. 5.4]{S}. See also
\cite[Corollary 2.7]{Z}.

\begin{cor} Let $K$ be an $H$-simple left $H$-comodule
algebra. If $W$, $U$ are left $K$-modules  then  $\dim K $ divides
$\dim W \dim U \dim H$. \qed
\end{cor}

Skryabin shows in \cite[Prop. 5.4]{S} that $\dim K \dim \End_K W$
divides $\dim W \dim U \dim H$ under the assumption $W$
irreducible, with a different proof.

\subsection{Examples of Yan-Zhu
stabilizers: Hopf subalgebras and left coideal
subalgebras}\label{subsection:leftcoidealsubalgebras} Let us
compute the Yan-Zhu stabilizers in some examples.

\begin{exa}\label{grupo} Let $G$ be a finite group and $H$ be the
group algebra of $G$. Let $F$ be a subgroup of $G$ and $\sigma\in
Z^2(F,\ku^{\times})$ be a normalized 2-cocycle. The twisted group
algebra $\ku_{\sigma} F$ is a left $H$-comodule algebra via
$\delta:\ku_{\sigma} F\to H\ot\, \ku_{\sigma} F$, $\delta(f)=f\ot
f$, for all $f\in F$.

For any left $\ku_{\sigma} F$-module $V$ the space $\End(V)$ is a
left $\ku F$-module via $ (f\cdot T)(v)=f\cdot T(f^{-1}\cdot v)$,
$T\in \End(V)$, $f\in F$, $v\in V$. The space $\ku G\ot_{\ku F}
\End(V)$ is a left $H$-module algebra and there is an isomorphism
of left module algebras
\begin{equation}\label{eqn:stab-group}
\St_{\ku_{\sigma} F}(V)\cong \ku G\ot_{\ku F} \End(V).
\end{equation}

\end{exa}

\begin{proof} If $g\in G, T\in\End(V)$ we denote by $\overline{g\ot T}$
the class of $g\ot T$ in $\ku G\ot_{\ku F} \End(V)$. Let
$\{x_i\}_{i\in I}$ be a complete set of representatives of the
right cosets. The $H$-module algebra structure in $\ku G\ot_{\ku
F} \End(V)$ is as follows; the action of $G$ is on the first
tensorand, and the product is given by
$$(\overline{x_i\ot\, T})(\overline{x_j\ot\, U})=\delta_{i,j}(\overline{x_i\ot\, T\circ U}), $$
for all $i,j\in I, T, U\in \End(V)$. We claim that the maps
$$\psi:\St_{\ku_{\sigma} F}(V)\to  \ku G\ot_{\ku F} \End(V),\quad \theta:\ku G\ot_{\ku F} \End(V)
\to \St_{\ku_{\sigma} F}(V)$$ defined by
$$\psi(\alpha_j\ot T_j)=\sum_{i} \alpha_j(x_i^{-1})\, \overline{x_i\ot T_j},
\quad \theta(\overline{g\ot T})=\sum_{f\in F}\delta_{fg^{-1}}\ot
f\cdot T,
$$
are well defined algebra isomorphisms, one the inverse of each
other. Indeed if $h\in F$ then
\begin{align*} \theta(\overline{gh\ot T})&= \sum_{f\in F}\delta_{fh^{-1}g^{-1}}\ot
f\cdot T=\sum_{f\in F}\delta_{fg^{-1}}\ot fh\cdot
T=\theta(\overline{g\ot h\cdot  T}).
\end{align*}
That is, $\theta$ is well-defined. Let $T, U\in \End(V)$ and
$\alpha_j\ot T_j, \beta_k\ot U_k \in \St_{\ku_{\sigma} F}(V)$ then
\begin{align*} \theta(x_i\ot\, T)\theta(x_j\ot\, U)&=\sum_{f, h\in F} (\delta_{fx_i^{-1}}\ot
f\cdot T )(\delta_{hx_j^{-1}}\ot h\cdot U)\\
&=\delta_{i,j}\sum_{f\in F}\delta_{fx_i^{-1}}\ot (f\cdot T)(f\cdot
U)=\theta((x_i\ot\, T)(x_j\ot\, U)),
\end{align*}

\begin{align*} \psi(\alpha_j\ot T_j)\psi(\beta_k\ot U_k)&= \sum_{i,l} \alpha_j(x_i^{-1})
 \beta_k(x_l^{-1})\, (\overline{x_i\ot T_j})(\overline{x_l\ot U_k})\\
 &=\sum_{i}\alpha_j(x_i^{-1})\beta_k(x_i^{-1})\, (\overline{x_i\ot T_j U_k})\\
 &=\psi((\alpha_j\ot T_j)(\beta_k\ot U_k)),
\end{align*}
thus $\theta$ and $\psi$ are algebra morphisms. Now let us compute
$\theta\circ \psi$ and  $ \psi\circ \theta$:

\begin{align*} \theta(\psi(\alpha_j\ot
T_j))&=\sum_{i}\alpha_j(x_i^{-1})\, \theta(\overline{x_i\ot T_j})
=\sum_{i, {f\in F}} \alpha_j(x_i^{-1})\, (\delta_{fx_i^{-1}}\ot
f\cdot T_j )\\
&=\sum_{i, {f\in F}} \alpha_j(f x_i^{-1})\, (\delta_{fx_i^{-1}}\ot
T_j )=\sum_{g\in G} \alpha_j(g)\, (\delta_g\ot T_j)=\alpha_j\ot
T_j.
\end{align*}
The third equality by \eqref{stab1}. On the other hand if $g=x_j
h, h\in F$ then
\begin{align*}\psi(\theta(\overline{g\ot T}))&=\sum_{f\in F}\psi(\delta_{fg^{-1}}\ot
f\cdot T)=\sum_{i, f\in F} \delta_{fg^{-1}}(x_i^{-1})\,
\overline{x_i\ot f\cdot T} &=\overline{x_j\ot h\cdot
T}=\overline{g\ot T}.
\end{align*}
\end{proof}

\begin{exa}\label{coideal-subalgebra} Let $H$ be a finite
dimensional Hopf algebra and $K\subseteq H$ be a left coideal
subalgebra, and therefore $K$ is a left $H$-comodule algebra via
the comultiplication. Denote $\overline{K}=H/\Ss^{-1}(K^{+})H$.
The canonical projection $\pi:H\to \overline{K}$ is an $H$-module
coalgebra map. The transpose of $\pi$ is an injective algebra
homomorphism $\overline{K}^{*}\hookrightarrow H^{*}$.  Via
$\pi^{*}$ the space $\overline{K}^{*}$ can be identified with the
subalgebra of $H^{*}$ consisting of elements $\alpha\in H^{*}$
such that $\alpha(x)=0$ for all $x\in \Ss^{-1}(K^{+})H$; clearly
this is a right coideal subalgebra of $H^{*}$. Let $V=\ku$ be the
trivial $K$-module. Then there is an isomorphism of right
$H^{*}$-comodule algebras
\begin{equation}\label{eqn:stab-coidealsubalg}\St_{K}(\ku)\cong \overline{K}^{*}
\end{equation}
\end{exa}
\begin{proof}  Since $\End(V)\simeq \ku$ we will identify
$\St_{K}(\ku)$ with a subalgebra of $H^{*}$. If
$\alpha\in\St_{K}(\ku)$ identity \eqref{stab1} implies that
$\varepsilon(k)<\alpha,t>= <\alpha,\Ss^{-1}(k)t>$, for any $t\in
H, k\in K$. Thus if $\varepsilon(k)=0$ then
$<\alpha,\Ss^{-1}(k)t>=0$, and therefore
$\alpha\in\overline{K}^{*}$. Reciprocally, if
$\alpha\in\overline{K}^{*} $ then
$\alpha(\Ss^{-1}(k)t-\varepsilon(k)t)=0$, since
$\Ss^{-1}(k)t-\varepsilon(k)t\in \Ss^{-1}(K^{+})H$,  and
\eqref{stab1} is fulfilled. This implies that
$\alpha\in\St_{K}(\ku)$.
\end{proof}

\subsection{Examples of Yan-Zhu
stabilizers: Yan-Zhu Stabilizers for Hopf Galois extensions} In
this subsection we shall give another expression for the Yan-Zhu
stabilizer in the case that $K$ is a Hopf-Galois extension over a
Hopf subalgebra $H'$ of $H$. First we recall the notion of
Hopf-Galois extensions.

Let $H'$ be a finite dimensional Hopf algebra.
\begin{defi} Let $K$ be a left $H'$-comodule algebra.
Set $R=K^{{\rm co} H'}$. The canonical map $\beta:K\ot_R K\to H\ot
K$ is defined  by $\beta(x\,\ot y)=x\_{-1}\ot x\_0y$, for all
$x,y\in K$. $K$ is called a {\it Hopf-Galois extension of $R$}
over $H'$ if $\beta$ is bijective.
\end{defi}

Following \cite{Sch} if $K\supseteq R$ is a Hopf-Galois extension
denote $\beta^{-1}(h\ot 1):=h^{[1]}\ot h^{[2]} \in K\ot_R K$, for
all $h\in H'$. The next result is due to H.-J. Schneider, see
\cite[Rmk 3.4]{Sch}.

\begin{lema} Let $K\supseteq R$ be a Hopf-Galois extension, then for
all $h,t\in H'$, $k\in K$, $r\in R$ we have that

\begin{align}\label{g1} r h^{[1]}\ot h^{[2]}&= h^{[1]}\ot
h^{[2]}r,\\
\label{g2} (th)^{[1]}\ot (th)^{[2]}&=t^{[1]}h^{[1]}\ot
h^{[2]}t^{[2]},\\
\label{g3} h^{[1]}h^{[2]}&=\varepsilon(h) 1_K,\\
\label{g4} h^{[1]}\ot 1\ot
h^{[2]}&=h\_1^{[1]}\ot h\_1^{[2]}h\_2^{[1]}\ot h\_2^{[2]}, \\
\label{g5} k\ot 1&= k\_{-1}^{[1]}\ot k\_{-1}^{[2]} k\_0,\\
\label{g6}  1\ot k& =k\_0 \Ss^{-1}(k\_{-1})^{[1]}\ot
\Ss^{-1}(k\_{-1})^{[2]},\\
\label{g7} h\_2\ot h\_1^{[1]}\ot
h\_1^{[2]}&=\Ss^{-1}(h^{[2]}\_{-1})\ot h^{[1]}\ot h^{[2]}\_0.
\end{align}
\end{lema}

\begin{proof} Equations \eqref{g1}, \eqref{g2}, \eqref{g5} and \eqref{g6}
follows by applying $\beta$. Since $(\varepsilon\ot\id)\beta=m$,
equation \eqref{g3} follows. To get \eqref{g4} apply
$(\id\ot\beta)(\beta\ot\id)$ on both sides. Finally, equation
\eqref{g7} follows from colinearity of $\beta$. More precisely,
$$\delta \beta=(\id_H\ot\beta)\widetilde{\delta}, $$
where $\delta:H\ot K\to  H\ot K\ot H$, $\widetilde{\delta}:K\ot_R
K\to K\ot_R K\ot H$ are defined by
$$\delta(h\ot y)=h\_1\ot y\_0\ot\, \Ss^{-1}(y\_{-1})h\_2, \quad\,
\widetilde{\delta}(x\,\ot y)=x\,\ot y\_0\ot\, \Ss^{-1}(y\_{-1}),$$
for all $x,y \in K, h\in H$.
\end{proof}

The following result is \cite[Corollary 3.5]{Sch}.

\begin{lema}\label{alg} Let $W$ be a left $K$-module. $\End_R(W)$ is a
left $H'$-module with respect to the action defined by
$$(h\cdot T)(w):= h^{[1]}T(h^{[2]}\,w), $$
for all $T\in\End_R(W)$, $h\in H'$, $w\in W$. We have also that
$\End_R(W)^{H'}=\End_K(W)$.
\end{lema}

\begin{proof} For any $h\in H'$, $T\in\End_R(W)$, $h\cdot T$ is an
$R$-module map by equation \eqref{g1}. Equation \eqref{g2} implies
that it is in fact an action. Equation \eqref{g3} and  \eqref{g4}
implies that $\End_K(W)$ is a module algebra over $H'$. Using
equation \eqref{g5} one can easily prove that the invariants of
$\End_R(W)$ are those who preserve the action of $K$.
\end{proof}

\smallbreak Let $H'$ be a Hopf subalgebra of $H$. Consider $H$ as
a left $H'$-module via the left regular representation. Let
$(x_j)_j\subseteq H$, $(\beta_j)_j\subseteq H^{*}$ be dual basis.

\begin{teo}\label{galois} Let $K\supseteq R$ be a Hopf-Galois
extension of $H'$. Let $W$ be a left $K$-module, then
\begin{enumerate}
   \item $\Hom_{H'}(H, \End_R(W))$ is a right $H^{*}$-comodule
   algebra, and
   \item $\St_K(W)\cong \Hom_{H'}(H, \End_R(W))$ as right $H^{*}$-comodule
   algebras.
\end{enumerate}
\end{teo}

\begin{proof} (1) The product is given by the convolution, that is if
$T,U\in\Hom_{H'}(H, \End_R(W))$ then $(TU)(h)=T(h\_1)U(h\_2)$, for
all $h\in H$. The identity is given by $\varepsilon$. The left
$H$-module structure $(h\cdot T)(x)=T(xh)$, $h, x\in H$,
$T\in\Hom_{H'}(H, \End_R(W))$ induces a right $H^{*}$-comodule
structure and becomes into a right $H^{*}$-comodule algebra.

(2) Let $\phi:\St_K(W)\to\Hom_{H'}(H, \End_R(W))$,
$\psi:\Hom_{H'}(H, \End_R(W))\to \St_K(W)$ be given by
$$\phi(\alpha_i\ot f_i)(h)(w)=\alpha_i(h) f_i(w), \quad\, \psi(
T)=\sum_j \beta_j\ot T(x_j), $$ for all $\alpha_i\ot f_i\in
\St_K(W)$, $h\in H, w\in W$. Let us verify that these maps are
well defined. Let $r\in R$, $t\in H'$, $h\in H$ and $w\in W$. Then
\begin{align*} \phi(\alpha_i\ot f_i)(h)(r\cdot w)&=\alpha_i(h)
f_i(r\cdot w)=\alpha_i(\Ss^{-1}(r\_{-1}h) r\_0\cdot f_i(w)
=\alpha_i(h)\, r\cdot f_i( w),
\end{align*}
the second equation by \eqref{stab1} and the last one because
$r\in R=K^{co H'}$. This proves that $\phi(\alpha_i\ot f_i)(h)$ is
an $R$-module map. We have also that

\begin{align*} t\cdot \phi(\alpha_i\ot f_i)(h)(w)&=t^{[1]}\phi(\alpha_i\ot f_i)
(h)(t^{[2]}\cdot w)=\alpha_i(h)\; t^{[1]}\cdot f_i(t^{[2]}\cdot w)\\
&=\alpha_i(\Ss^{-1}(t^{[2]}\_{-1})h)\; t^{[1]}t^{[2]}\_0\cdot f_i(
w)\\
&=\alpha_i(t\_2h)\; t\_1^{[1]}t\_1^{[2]}\cdot f_i(
w)=\alpha_i(th)f_i( w).
\end{align*}
The third equation by \eqref{stab1}, the fourth by \eqref{g7} and
the fifth by \eqref{g3}. This proves that $\phi(\alpha_i\ot f_i)$
is an $H'$-module map and therefore $\phi$ is well defined. The
proof that $\psi(T)\in \St_K(W)$ is done using \eqref{g6} and
\eqref{stab1}. That $\phi$ is an algebra map and a right
$H^{*}$-comodule morphism is a straightforward computation. The
identities $\psi\phi=\id, \phi\psi=\id$ are checked without
difficulties. \end{proof}

\section{Applications of the Yan-Zhu stabilizers to module categories}\label{section:modcat}

\subsection{Internal Hom}\label{subsection:ends}
We keep the notation of the preceding section.

\begin{prop}\label{stab=internalend} $\underline{\Hom} (U, W) =\St_K(U, W) $, and the bilinear map
$$\underline{\Hom} (U, W) \times \underline{\Hom} (W, Z) \to \underline{\Hom} (U, Z)$$
coincides with \eqref{stabcomp}.
\end{prop}

\pf Let us identify $H^* \ot \Hom (U, W)$ with $\Hom(H\ot U, W)$ in the natural way. Let $X\in \Rep H$. There are natural linear inverse isomorphisms
\begin{align*}G: \Hom_H(X, H^* \ot \Hom (U, W)) &\to \Hom(X\otimes U, W),
\\ F:  \Hom(X\otimes U, W) &\to \Hom_H(X, H^* \ot \Hom (U, W)),
\end{align*}
given by $G(\psi) (x\,\ot u) = \psi(x) (1\ot u)$,  $F(\phi)(x)
(a\ot u)  = \phi(a\cdot x\,\ot u)$. Here and in the rest of the
proof, $\psi \in \Hom_H(X, H^* \ot \Hom (U, W))$, $\phi \in
\Hom(X\otimes U, W)$, $x\in X$, $a\in H$, $u\in U$; and also
$\alpha\in H^*$, $k\in K$. Given $\psi$ and $x$, we write
symbolically $\psi(x) = \psi(x)\_1 \ot \psi(x)\_2$, with
$\psi(x)\_1 \in \End (H^*)$, $\psi(x)\_2 \in \Hom(U,W)$. Note that
$$
\ele(\psi(x)) (\alpha \otimes u) (a) = \langle \alpha, a\_2\rangle \psi(x)(a\_1 \ot u).
$$

We claim that  $F(\Hom_K(X\otimes U, W)) = \Hom_H(X, \St_K(U,
W))$, up to identification via $\ele$. Indeed, $\phi \in
\Hom_K(X\otimes U, W)$ iff $\psi = G(\phi)$ satisfies
\begin{equation}\label{phiequiv}
k\cdot \left(\psi(x)(a\ot u)\right) = \psi(x)(k\_{-1}a\ot k\_0\cdot u).
\end{equation}

\smallbreak Let us denote by $\rho$ either $\rho_{H^*\ot W}: K\to
\End H^* \ot \End W \simeq \End(H^* \ot W)$ or $\rho_{H^*\ot U}$.
We compute on one hand
\begin{align*}
\left[\psi(x) \rho(k)(\alpha \ot u)\right](a)
&= \left[\psi(x) (k\_{-1}\rightharpoondown\alpha \ot k\_0\cdot u)\right](a)
\\
  &= \langle k\_{-1}\rightharpoondown\alpha, a\_2\rangle \psi(x)(a\_1 \ot k\_0\cdot u)
\\
  &= \langle \alpha, \Ss^{-1}(k\_{-1})a\_2\rangle \psi(x)(a\_1 \ot k\_0\cdot u),
\\
  &= \boxtimes_1.
\end{align*}
and on the other

\begin{align*}
\left[\rho(k)\psi(x) (\alpha \ot u)\right](a)
&=  \left[k\_{-1}\rightharpoondown\psi(x)\_1(\alpha) \ot k\_0\cdot \psi(x)\_2(u)\right](a)
\\
  &=  \left[\psi(x)\_1(\alpha) \ot k\_0\cdot \psi(x)\_2(u)\right](\Ss^{-1}(k\_{-1})a)
\\
  &= \langle \alpha, \Ss^{-1}(k\_{-2})a\_{2}\rangle \, k\_0\cdot
\left[\psi(x)(\Ss^{-1}(k\_{-1})a\_{1} \ot u)\right]
\\
  &= \boxtimes_2.\end{align*}

If $\boxtimes_1 = \boxtimes_2$ then taking $\alpha = \varepsilon$ we
get \eqref{phiequiv}. Conversely, it is not difficult to see that
\eqref{phiequiv} implies $\boxtimes_1 = \boxtimes_2$.

\smallbreak Now the claim says that  $\Hom_K(X\otimes U, W) \simeq
\Hom_H(X, \St_K(U,W))$; so that the functor $X\mapsto
\Hom_K(X\otimes U, W)$ is representable by $\St_K(U,W)$.
Furthermore, it is not difficult to see that the composition
\eqref{stabcomp} satisfies the defining property in \cite[Section
3.3]{O1}, and the proposition follows. \epf

\subsection{Exact module categories}
We state our first application of Proposition
\ref{stab=internalend}.

\begin{cor}\label{exact=lca} Assume that $K$ is exact. Let $W$ be a generator
of $_K\Mo$. Then ${}_K\Mo$ is equivalent to $_{H}\Mo_{\St_K(W)}$
as module categories over $\Rep H$.
\end{cor}

\pf This follows from Theorem \ref{exactindecompmc} by Proposition
\ref{stab=internalend}. \epf

We now give a refinement of Proposition \ref{todas}.

\begin{teo}\label{mod=lca} Any indecomposable exact module category over
$\Rep H$ is equivalent to ${}_K\Mo$ for some $H$-simple from the
right left $H$-comodule algebra $K$, with $K^{\co H}\simeq \ku$.
\end{teo}

\pf By Theorem \ref{exactindecompmc}, there exists a $H$-module
algebra $R$ such that
\begin{equation}\label{one}
\Mo\simeq\, _{H}\Mo_R
\end{equation}
as module categories over $\Rep H$. Because of \cite[Lemma
4.2]{eo}, see Remark \ref{irredimpliesHsimple}, we can assume that
$R$ has no non trivial $H$-stable ideals. Hence ${}_{R}\Mo$ is
exact as module category over $\Rep (H^*)^{\cop}$ by Proposition
\ref{Hsimplealgexact}. It follows that $R$ is quasi-Frobenius by
Remark \ref{exactmc-k-frob}.

Let $W$ be a simple $R$-module and set $K=\St_R(W)$. We know that
$W$ generates $_{H}\Mo_R$, see the remarks previous to Theorem
\ref{exactindecompmc}. Also, $K$ is $H$-simple from the right
because of Remark \ref{irredimpliesHsimple}. Hence ${}_{K}\Mo$ is
indecomposable by Proposition \ref{decomposable}, and exact by
Proposition \ref{Hsimplealgexact}. Again, $W\neq 0$ is a generator
of ${}_{K}\Mo$. Observe next that
\begin{equation}\label{two}
{}_{K}\Mo \simeq\, _{H}\Mo_{\St_K(W)}
\end{equation}
by Corollary \ref{exact=lca}. Now  $R\simeq \St_{\St_R(W)}(W)$ by
Theorem \ref{cor:yzduality-simple} (Yan-Zhu duality). The
Proposition now follows from this, \eqref{one} and \eqref{two}.
\epf

\subsection{The dual module category}
Another important tool in the study of tensor categories is the
notion of \emph{dual tensor category} with respect to a module
category. In some sense this notion is the categorification of the
notion of a centralizer of an algebra. The dual tensor category
has been intensively used in \cite{ENO}. See also \cite{O1},
\cite{O2}.

\smallbreak Let $\ca$ be a finite tensor category.

\begin{defi} Let $\Mo$ be an exact module category over $\ca$. The
\emph{dual tensor category}  (with respect to $\Mo$) is the
category $\ca^*_{\Mo} := \Fun_{\ca}(\Mo,\Mo)$ with the tensor
product given by the composition of module functors
\eqref{modfunctor-comp}.
\end{defi}

If $\No$ is a module category over $\ca$ then
$\Fun_{\ca}(\No,\Mo)$ is a module category over $\ca^*_{\Mo}$ via
the composition $\Fun_{\ca}(\Mo,\Mo)\times \Fun_{\ca}(\No,\Mo)\to
\Fun_{\ca}(\No,\Mo)$, see \eqref{modfunctor-comp} again.

\begin{prop}\label{corr} Let $\No$ be an exact module
category over $\ca$. Then
\begin{enumerate}
    \item $\Fun_{\ca}(\No,\Mo)$ is an exact module category over
    $\ca^*_{\Mo}$.
    \item The map  $\No \mapsto \Fun_{\ca}(\No,\Mo)$ is
    a bijective correspondence between equivalence classes of
    exact module categories over $\ca$ and $\ca^*_{\Mo}$.
\end{enumerate}
\end{prop}
\pf See \cite[Lemma 3.30]{eo} and \cite[Theorem 3.31]{eo}. \epf

\begin{lema}\label{corr2} Let $A\in\ca$ be an algebra and assume that $\Mo=\ca_A$ is an exact
module category over $\ca$. Then
\begin{enumerate}
    \item The tensor categories $\ca^*_{\Mo}$ and $({}_A\ca_A)^{\op}$ are
    equivalent.
    \item The bijective correspondence between equivalence classes of
    exact module categories over $\ca$ and $({}_A\ca_A)^{\op}$ arising from
    Proposition \ref{corr} (2) is explicitly given by
    $\ca_B \mapsto {}_B\ca_A$, $B$ any algebra in $\ca$.
\end{enumerate}
\end{lema}
\pf See the proof of \cite[Lemma 3.30]{eo}.\epf

\begin{rmk}\label{mod-op} If $A = \uno$, then we conclude that
the tensor categories $\ca^*_{\ca}$, $\ca^{\op}$ are equivalent.
Hence the correspondence of exact module categories over $\ca$ and
$\ca^{\op}$ is just $\ca_B\mapsto {}_B\ca$, $B$ an algebra in
$\ca$.
\end{rmk}

\subsection{ Correspondence of module categories over $\Rep(H^*)$
and $\Rep H $} Let $H$ be a finite-dimensional Hopf algebra. In
this subsection we study a bijective correspondence between
equivalence classes of exact module categories over $\Rep H $ and
$\Rep(H^*)$, and show that this agrees with Proposition \ref{corr}
(2). Roughly the correspondence is as follows.  If $K$ is a
$H$-simple left $H$-comodule algebra, then $K^{\op}$ is a left
$H^*$-module algebra and therefore is a right $H$-comodule
algebra. If $V\neq 0$ is a right $K$-module, then the stabilizer
$\St_{K^{\op}}(V)$ is a left $H^*$-comodule algebra and the module
category ${}_{\St_{K^{\op}}(V)}\Mo$ does not depend on $V$.
Therefore we have a map
$${}_K\Mo\longmapsto\;
{}_{\St_{K^{\op}}(V)}\Mo$$ assigning module categories over $\Rep
H $ to module categories over $\Rep(H^*)$.

\smallbreak We begin by the following well-known Lemma. Recall
that $H^*$ is a left $H$-module algebra via $\rightharpoonup$.

\begin{lema}\label{equivalence} There is a tensor equivalence
$\Rep(H^*)\simeq  {}_{H^*\!}\Rep H _{H^*}$.
\end{lema}
\pf We only sketch the proof. The functors $\Fc:\Rep(H^*)\to
{}_{H^*\!}\Rep H _{H^*}$, \newline $\Gc:{}_{H^*\!}\Rep H _{H^*}\to
\Rep(H^*)$ are defined by $\Gc(V)= V^{H}=\{v\in V: h\cdot v=
\varepsilon(h) v\}$, $\Fc(X)=X\ot H^*$ with the following
structure. For all $x\in X, \alpha, \beta\in H^*, h\in H$
\begin{align*} h\cdot (x\ot \alpha)=x\ot h\rightharpoonup
\alpha,\quad \beta\cdot (x\ot \alpha)= \beta\_1 \cdot x\ot
\beta\_2 \alpha,\quad (x\ot \alpha)\cdot \beta=x\ot \alpha\beta.
\end{align*}
\epf

The next Proposition shows that there is a correspondence between
module categories over $\Rep H $ and $\Rep(H^*)$. This result was
first established by Ostrik for weak Hopf algebras. See
\cite[Theorem 5]{O2}.

\begin{prop}\label{corr3} There is a bijective correspondence between
equivalence classes of module categories over $\Rep H $ and
$\Rep(H^*)$. Explicitly, if $R$ is  $H$-module algebra, then $\Rep
H _R\mapsto {}_{H^*}\Rep H _R$, where $\ot: \Rep (H^*) \times
{}_{H^*}\Rep H _R \to {}_{H^*}\Rep H _R$ is given by $X\ot V =
\Fc(X) \ot_{H^*} V$.
\end{prop}

\pf We know that the module category $\Vect \simeq \Rep H_{H^*}$.
By Lemma \ref{corr2} applied to $A = H^*$, we have $(\Rep
H)^*_{\Vect} \simeq ({}_{H^*\!}\Rep H _{H^*})^{\op}$. Hence the
application we are looking into is just the composition
$$
\xymatrix{\mod \Rep H \ar[dr]_{\text{Remark \ref{mod-op}}\qquad
}\ar[rr] && \mod
({}_{H^*\!}\Rep H _{H^*}) \\
&\mod (\Rep H)^{\op}   \ar[ur]_{\qquad \text{Lemma \ref{corr2} (2)}}& \\
 \Rep H_R \ar@{|->}[rr]\ar@{|->}[dr] && {}_{H^*\!}\Rep H _{R}
\\&{}_R\Rep H \ar@{|->}[ur] & }
$$
combined with Lemma \ref{equivalence}. \epf

The main result of this Subsection is an explicit description of
this correspondence.

\begin{teo}\label{corr4} The bijective correspondence between equivalence classes of
exact module categories over $\Rep H $ and $\Rep(H^*)$ settled in
Proposition \ref{corr3} coincides with the map
$${}_K\Mo\longmapsto\;
{}_{\St_{K^{\op}}(V)}\Mo,$$ $K$ an $H$-simple left $H$-comodule
algebra and $V\neq 0$ a right $K$-module.
\end{teo}

In presence of Example \ref{coideal-subalgebra}, the Theorem
``explains" the correspondence between coideal subalgebras of $H$
and $H^*$ described by Masuoka \cite[2.10 (iii)]{Mk}.

\pf Let $K$ be a $H$-simple left $H$-comodule algebra and let $W$
be a left $K$-module. By Corollary \ref{exact=lca} there is an
equivalence of module categories $(F,c): {}_K\Mo\to \Rep H
_{\St_K(W)}$.

We first claim that the functor $F$ induces an equivalence
\begin{equation}\label{corr4-eqn}
\Rep(H^*)_{K^{\op}}\simeq {}_{H^*}\Rep H _{\St_K(W)}
\end{equation} of module
categories over $\Rep(H^*)$. We first observe that indeed
$K^{\op}$ is a left $H^*$-module algebra with left action given by
$\alpha\cdot  x= \langle \alpha, \Ss^{-1}(x\_{-1})\rangle\, x\_0$,
for every $\alpha \in H^*$, $x\in K$. Thus, we may consider the
module category $\Rep(H^*)_{K^{\op}}$ over $\Rep(H^*)$. Now, an
object $M\in\Rep(H^*)_{K^{\op}}$ is a left $K$-module $\cdot:K\ot
M\to M$ provided with a left $H^*$-action $\triangleright :H^*\ot
M\to M$ such that
\begin{equation}\label{equiv-formula1} \alpha\triangleright
(x\cdot m)=\langle\alpha\_2, \Ss^{-1}(x\_{-1})\rangle\,
 x\_0\cdot (\alpha\_1\triangleright m)
\end{equation}
Recall that $H^*$ is a left $H$-module via $\rightharpoonup$ and
if $M$ is a $K$-module then $H^*\ot M\in {}_{K}\Mo$.
\begin{step}\label{equiv-1}
The left action $\triangleright:H^*\ot M\to M$ is a $K$-module
map.
\end{step}
\pf The map $\triangleright:H^*\ot M\to M$ is a $K$-module map if
and only if
\begin{equation}\label{equiv-formula2} k\cdot
(\alpha\triangleright m) = (k\_{-1}\rightharpoonup
\alpha)\triangleright (k\_0 \cdot m),
\end{equation}
for all $\alpha\in H^*, m\in M, k\in K$. The right hand side
equals to
\begin{align*} &=\langle\alpha\_2, k\_{-1}\rangle\; \alpha\_1\triangleright (k\_0 \cdot
m)\\
&=\langle\alpha\_3, k\_{-2}\rangle\langle\alpha\_2,
\Ss^{-1}(k\_{-1})\rangle\; k\_0\cdot (\alpha\_1\triangleright m)\\
&=\langle\alpha\_2,\Ss^{-1}(k\_{-1}) k\_{-2}\rangle\;k\_0\cdot (\alpha\_1\triangleright m)\\
&=k\cdot (\alpha\triangleright m).
\end{align*}
The second equality follows from \eqref{equiv-formula1}.\epf

\smallbreak Step \ref{equiv-1} says that the map
$\triangleright:H^*\ot M\to M$ is in the category ${}_K\Mo$. Thus,
applying the functor $F$, we get the map $F(H^*\ot M)
\xrightarrow{\,F(\triangleright)\,} F(M)$; we can consider the
composition
\begin{equation}\label{dual-action}H^*\ot F(M)\xrightarrow{\,\,\,c^{-1}_{H^*,M}\,}
F(H^*\ot M)\xrightarrow{\,\,F(\triangleright)\,}
F(M).\end{equation}

\begin{step}\label{equiv-2} Suppose $M\in \Rep(H^*)_{K^{\op}}$. Then \begin{enumerate}
    \item The composition \eqref{dual-action} is a left $H^*$-action on
    $F(M)$.
    \item $F(M)$ is an object in $ {}_{H^*}\Rep H _{\St_K(W)}$.
\end{enumerate}
\end{step}
\pf (2) follows from (1), since the composition
\eqref{dual-action} is a morphism in  $\Rep H _{\St_K(W)}$. Let us
prove (1). The associativity of the action given by
\eqref{dual-action} is equivalent to
$$F(\triangleright)\,
c^{-1}_{H^*,M}\, (\mu\ot \id_{F(M)})= F(\triangleright)\,
c^{-1}_{H^*,M}\left(\id_{H^*}\ot F(\triangleright)\,
c^{-1}_{H^*,M}\right),$$ where $\mu$ denotes the multiplication of
$H^*$. Since $\mu:H^*\ot H^*\to H^*$ is a morphism in $\Rep H $
the naturality of $c$ implies that
\begin{equation}\label{t11}(\mu\ot \id_{F(M)})\, c_{H^*\ot H^*,M}=
c_{H^*,M}\, F(\mu\ot \id_M).\end{equation} Analogously, since
$\triangleright:H^*\ot M\to M$ is in ${}_K\Mo$ the naturality of
$c$ implies that
\begin{equation}\label{t12}(\id_{H^*}\ot F(\triangleright))\,
 c_{H^*, H^*\ot M}= c_{H^*, M}\, F(\id_{H^*}
\ot\, \triangleright).
\end{equation}
Hence
\begin{align*} F(\triangleright)\,
c^{-1}_{H^*,M}\, (\mu\ot \id_{F(M)})&=F(\triangleright)F(\mu\ot
\id_M)\, c^{-1}_{H^*\ot H^*,M}\\
&=F(\triangleright\,(\mu\ot \id_M))\,  c^{-1}_{H^*, H^*\ot
M}(\id_{H^*}\ot \,c^{-1}_{H^*,M})\\
&=F(\triangleright) F(\id_{H^*}\ot \, \triangleright)\,
c^{-1}_{H^*, H^*\ot
M}(\id_{H^*}\ot \,c^{-1}_{H^*,M})\\
&=F(\triangleright)\, c^{-1}_{H^*,M} (\id_{H^*}\ot
F(\triangleright)) (\id_{H^*}\ot \,c^{-1}_{H^*,M})
\end{align*}
The first equality by \eqref{t11}, the second by
\eqref{modfunctor1}, the third because $\triangleright$ is an
action and the last equality follows from \eqref{t12}. \epf

\begin{step}\label{morfismos} Let $X\in \Rep(H^*)$, $M\in
\Rep(H^*)_{K^{\op}}$. The map $c_{X,M}: F(X\ot M)\to X\ot F(M)$ is
a morphism in  ${}_{H^*}\Rep H _{\St_K(W)}$. Here $X$ is an
$H$-module with trivial action.
\end{step}
\pf By definition the map $c_{X,M}$ is a morphism in $\Rep H
_{\St_K(W)}$. Thus, we only must show that $c_{X,M}$ is a morphism
of $H^*$-modules. This is equivalent to prove that
\begin{multline}\label{morfismos1} c_{X,M}
F(\theta)c^{-1}_{H^*,X\ot M} \\ = (\triangleright_X\ot
F(\triangleright_M) c^{-1}_{H^*,M})(\id_{H^*}\ot
\tau\ot\id_{F(M)})(\Delta\ot\id_{X\ot F(M)})(\id_{H^*}\ot
c_{X,M}).
\end{multline}
Here,  $\triangleright_X$ and $\triangleright_M$ are the
$H^*$-actions on $X$, $M$ respectively; $\tau :H^*\ot X\to X\ot
H^*$ is the usual transposition; and $\theta:H^*\ot X\ot M\to X\ot
M$ is the left $H^*$-action on $X\ot M$, that is
$\theta=(\triangleright_X\ot\,\triangleright_M)(\id_{H^*}\ot
\tau\ot \id_M)(\Delta\ot \id_X\ot \id_M)$. Let $\phi:H^*\ot X\to
X\ot H^*$ be the morphism of $H$-modules defined by  $
\phi=(\id_{H^*}\ot\,
\triangleright_X)(\id_{H^*}\ot\tau_{H^*X})(\Delta\ot\id_X)$. By
the naturality of $c$ implies that the diagram
\begin{align}\label{diagram11}
\begin{CD}
 F(H^*\ot X\ot M) @> \text{id} \ot c_{H^*\ot X,M} >> H^*\ot X\ot F(M)\\
@V F(\phi\ot\text{id}) VV    @VV\phi\ot\text{id} V   \\
F(X\ot H^*\ot M) @>>c_{X\ot H^*, M}> X\ot H^*\ot F(M)
\end{CD}
\end{align}
is commutative. Since the map $\triangleright_M:H^*\ot M\to M$ is
a $K$-module map-- Step \ref{equiv-1}-- the naturality of $c$
implies the commutativity of the diagram
\begin{align}\label{diagram12}
\begin{CD}
 F(X\ot H^*\ot M) @>  c_{ X,H^*\ot M} >> X\ot F(H^*\ot M)\\
@V F(\text{id} \ot \triangleright_M) VV    @VV\text{id} \ot F(\triangleright_M) V   \\
F(X\ot  M) @>>c_{X*, M}> X\ot F(M).
\end{CD}
\end{align}
Then, \eqref{morfismos1} equals to
\begin{align*} &= (\id_X\ot F(\triangleright_M)
c^{-1}_{H^*,M})(\phi\ot \id) (\id\ot  c_{X,M})\\
&=(\id_X\ot F(\triangleright_M) c^{-1}_{H^*,M})(\phi\ot \id)
c_{H^*\ot X, M} c^{-1}_{H^*, X\ot M}\\
&=(\id_X\ot F(\triangleright_M) c^{-1}_{H^*,M})c_{X\ot H^*, M}
F(\phi\ot\id) c^{-1}_{H^*, X\ot M}\\
&=(\id_X\ot F(\triangleright_M))c_{X,H^*\ot M} F(\phi\ot\id) c^{-1}_{H^*, X\ot M}\\
&=c_{X*, M} F(\id\ot \triangleright_M) F(\phi\ot\id) c^{-1}_{H^*,
X\ot M}.
\end{align*}
The second and the fourth equalities by \eqref{modfunctor1}, the
third by \eqref{diagram11} and the last by \eqref{diagram12}. \epf

By Step \ref{equiv-2} (2), the restriction $F:\Rep(H^*)_{K^{\op}}
\to {}_{H^*}\Rep H _{\St_K(W)}$ is well-defined. If $X\in
\Rep(H^*)$, $M\in \Rep(H^*)_{K^{\op}}$, then the map
$d_{X,M}:F(X\ot M)\to (X\ot H^{*})\ot_{H^{*}} F(M)$ is defined as
the composition $F(X\ot M)\xrightarrow{\,c_{X,M} \,}X\ot
F(M)\xrightarrow{\,\simeq\,} (X\ot H^{*})\ot_{H^{*}} F(M)$. Here
$X$ is considered as a trivial left $H$-module. Step
\ref{morfismos} implies that $d_{X,M}$ is a morphism in
${}_{H^*}\Rep H _{\St_K(W)}$. This finishes the proof of
\eqref{corr4-eqn}. The Theorem now follows from the commutativity
of the following diagram:
$$
\xymatrix{{}_K\Mo \ar[rr]^{\sim\qquad}_{\text{Cor.
\ref{exact=lca}}\qquad}\ar[rrd]_{\sim\qquad} && \Rep
H_{\St_{K}(W)} \ar[rr]^{\sim\qquad}_{\text{Prop.
\ref{corr3}}\qquad} &&{}_{H^*}\Rep H _{\St_{K}(W)}\ar[d]^{\wr}_{
\eqref{corr4-eqn}}
\\ &&{}_{\St_{K^{\op}}(V)}\Mo && \ar[ll]_{\sim\qquad}^{\text{Cor.
\ref{exact=lca}}\qquad} \Rep(H^*)_{K^{\op}}}
$$
 \epf

\begin{cor} Let $K$ be a $H$-simple left $H$-comodule algebra. If
$V$ and $W$ are left $K$-modules then $\St_K(V)$ and $\St_K(W)$
are equivariant Morita equivalent. \qed \end{cor}

\end{document}